# Variational (Energy-Based) Spectral Learning: A Machine Learning Framework for Solving Partial Differential Equations

M. M. Hammad

Faculty of Science, Department of Mathematics and Computer Science, Damanhour University, Egypt
Email: m_hammad@sci.dmu.edu.eg
Orcid: https://orcid.org/0000-0003-0306-9719

## Abstract

We introduce variational spectral learning (VSL), a machine learning framework for solving partial differential equations (PDEs) that operates directly in coefficient space of spectral expansions. VSL offers a principled bridge between variational PDE theory, spectral discretization, and contemporary machine learning practice. The core idea is to recast a given PDE $\mathcal{L}u = f$ in $Q = \Omega \times (0, T)$ with boundary and initial conditions, into differentiable space–time energies built from strong-form least-squares residuals and weak (Galerkin) formulations. The solution is represented as a finite spectral expansion $u_N(x,t) = \sum_{n=1}^{N} c_n \varphi_n(x,t)$, where $\varphi_n$ are tensor-product Chebyshev bases in space and time, with Dirichlet-satisfying spatial modes enforcing homogeneous boundary conditions analytically. This yields a compact linear parameterization in the coefficient vector **c**, while all PDE complexity is absorbed into the variational energy. We show how to construct strong-form and weak-form space–time functionals, augment them with initial-condition and Tikhonov regularization terms, and minimize the resulting objective with gradient-based optimization. In practice, VSL is implemented in TensorFlow using automatic differentiation and Keras cosine–decay-with-restarts learning-rate schedules, enabling robust optimization of moderately sized coefficient vectors. Numerical experiments on benchmark elliptic and parabolic problems, including one- and two-dimensional Poisson, diffusion, and Burgers-type equations, demonstrate that VSL attains accuracy comparable to classical spectral collocation with Crank–Nicolson time stepping, while providing a differentiable objective suitable for modern optimization tooling. Comparisons with physics-informed neural networks highlight that learning in coefficient space mitigates some of the optimization pathologies of deep field-space surrogates, particularly for problems with smooth solutions and stiff operators.

## 1. Introduction

Partial differential equations (PDEs) are central to the mathematical formulation of phenomena in physics, engineering, finance, and many areas of applied science [1, 2]. The numerical solution of PDEs has traditionally relied on well-established discretization techniques such as finite difference, finite volume, finite element, and spectral methods, which have achieved a high level of maturity in terms of accuracy, robustness, and theoretical understanding [3–15]. In parallel, recent developments in machine learning (ML) have led to new paradigms for solving PDEs by embedding physical laws directly into learning objectives, exemplified by approaches such as physics-informed neural networks (PINNs) and operator-learning architectures [16–23]. These methods seek to leverage automatic differentiation, expressive function approximators, and large-scale optimization to construct data-efficient surrogates or direct solvers for PDEs [24–26]. Despite their promise, there remains a significant gap between classical variational formulations of PDEs and many of the current machine learning-based solvers, both in how the PDE constraints are enforced and in how approximation spaces are parameterized and optimized.

Among classical methods, variational formulations occupy a central position, particularly in elliptic and parabolic problems [27–30]. By recasting PDEs as minimization problems for energy functionals or as weak formulations in appropriate function spaces, they provide a rigorous foundation for Galerkin and Petrov–Galerkin methods, including spectral and finite element schemes [7–9], [11–15], [27–30]. These methods achieve high accuracy and stability by aligning the approximation spaces with the underlying functional analysis of the PDE, and by enforcing residual orthogonality in a weak sense. Spectral methods, in particular, exploit global polynomial or trigonometric bases to deliver exponential convergence for smooth solutions [11–15]. However, such methods typically rely on explicitly derived variational forms, hand-crafted bases, and direct or iterative solvers for linear or mildly nonlinear algebraic systems [27–30]. Their extension to complex, heterogeneous, or parameterized problems can be technically demanding and computationally costly, especially when repeated solves are required for many parameter configurations [31, 32].

On the other hand, learning-based PDE solvers often adopt a field-space viewpoint, in which a neural network parameterizes the solution $u(x)$ or $u(x,t)$ directly as a function of space (and time) [16–19]. The PDE is then incorporated through a residual loss, either in strong form (pointwise residuals of the differential operator) or in a weak-form integral sense, and the network parameters are tuned by gradient-based optimization [33, 34]. While such approaches benefit from automatic differentiation and flexible architectures, they face several challenges. Training dynamics can be stiff and highly sensitive to hyperparameters, particularly when high-order derivatives are involved or when the PDE exhibits multiple scales [35–39]. Boundary and initial conditions are frequently enforced by penalty terms rather than analytically, which may lead to constraint violations or the need for large penalty weights [40, 41]. Moreover, the underlying variational structure of the PDE is often only partially exploited: the choice of loss function does not always correspond to a principled energy functional, and convergence guarantees are typically weaker than those available for classical Galerkin schemes [33, 34], [42, 43].

These observations motivate a more systematic integration of variational PDE structure, spectral approximation, and modern machine learning tools. The central idea of this work is to return to the variational roots of PDE discretization while adopting the optimization and automatic differentiation machinery of contemporary deep learning. Instead of parameterizing the solution field directly by a neural network and enforcing the PDE through loosely variational losses, we propose to parameterize the solution in a spectral basis and define loss functions that arise directly from strong and weak variational formulations. The resulting framework, which we term variational spectral learning (VSL), can be viewed as gradient-based learning in coefficient space, guided by PDE-derived energies, and implemented with modern optimization algorithms and automatic differentiation.

In the VSL framework, the unknown solution $u(x, y, t)$ is expanded in a tensor-product spectral basis in space and time. For problems on simple domains such as rectangles, we employ Chebyshev polynomials in each spatial direction and in time, with spatial basis functions constructed to satisfy homogeneous Dirichlet boundary conditions analytically. For example, in the two-dimensional case, separate Dirichlet-satisfying bases in $x$ and $y$ are combined with a temporal Chebyshev basis to form a space–time spectral ansatz. The degrees of freedom of the approximation are the spectral coefficients of this expansion, which are treated as trainable parameters. This coefficient-space parameterization retains the interpretability and convergence properties of spectral methods while enabling the use of gradient-based optimization and generic autodifferentiation for assembling residuals and energies.

The variational character of VSL enters through the definition of training objectives. For a given PDE operator $\mathcal{L}$ and forcing term $f$, we consider two primary forms of variational energy. The strong-form energy is defined as a least-squares functional of the pointwise PDE residual [44, 45]. Writing the PDE abstractly as $F(u) = 0$ on the space–time domain, the strong energy integrates $|F(u)|^2$ over space–time, yielding a scalar objective closely related to strong-form least-squares discretizations. The weak-form energy is based on a Galerkin viewpoint. Depending on the PDE class and desired derivative order, it can be implemented either (i) as a classical weak form obtained by integration by parts (thereby reducing derivative order and aligning with coercive bilinear forms), or (ii) as a Galerkin/moment residual obtained by projecting the strong residual onto the spectral trial basis [11–15], [27–30]. In either case, the resulting residual quantities are assembled by quadrature and combined in a least-squares objective over the corresponding coefficient-wise residuals [46]. Both energies are approximated by tensor-product quadrature rules in space–time, and their gradients with respect to spectral coefficients are obtained automatically via nested differentiation tapes.

Boundary and initial conditions are treated in a manner consistent with the variational formulation. Homogeneous Dirichlet boundary conditions in space are built into the basis construction, so that any coefficient vector automatically yields a solution satisfying $u = 0$ on the boundary [11–15]. This eliminates the need for boundary penalty terms or explicit constraint enforcement in the optimization [40, 41]. Initial conditions in time-dependent problems are incorporated through additional loss terms that penalize deviations between the spectral approximation at $t = 0$ and a prescribed initial profile. The overall training objective therefore combines a variational PDE energy, an initial-condition loss, and a regularization term on the coefficients, yielding a well-structured optimization problem whose minima correspond to approximate solutions of the original PDE with correctly enforced constraints [27–30].

From an optimization standpoint, VSL recasts variational PDE solving as training a shallow, linear model in coefficient space under a physically motivated loss. This perspective invites the use of advanced optimization techniques developed in the machine learning community, such as adaptive gradient methods and learning-rate schedules with restarts [47–50]. In particular, this work employs a cosine-decay-with-restarts learning-rate schedule (e.g., as implemented in Keras), which modulates the learning rate in oscillatory cycles, providing a practical mechanism to modulate step sizes over training, which can mitigate stagnation during early iterations while supporting fine-scale refinement in later phases [49, 50]. Because the parameterization is linear in the coefficients, the map $\mathbf{c} \to u_N(\cdot; \mathbf{c})$ is a fixed linear feature map, where $u_N$ is finite spectral expansion and $\mathbf{c}$ is the coefficient vector, and for linear PDEs the resulting discrete objectives are typically convex quadratics under standard assumptions and sufficiently accurate quadrature [46], [51]. For nonlinear PDEs, nonconvexity remains (through the operator), but the coefficient-space formulation avoids the additional architectural nonlinearity introduced by deep field-space neural networks and is therefore often empirically easier to optimize. At the same time, VSL inherits the flexibility of autodifferentiation and can readily accommodate higher-order derivatives and complex PDE operators without manual derivation of discrete stencils; moreover, both strong residuals and Galerkin-type residual constructions can be encoded at the continuous level and differentiated automatically after quadrature discretization.

To assess the practical capabilities of VSL, it is essential to compare it against both classical and learning-based baselines. Classical spectral methods provide a natural benchmark because they share the same functional approximation spaces but rely on direct or implicit linear solvers rather than gradient-based minimization of a variational energy. For time-dependent problems, we consider Chebyshev–Lobatto collocation in space combined with Crank–Nicolson time stepping, a high-order, A-stable method that is widely used for diffusion-type equations [13], [52–55]. This solver discretizes the PDE into a linear system at each time step and enforces the boundary conditions directly at the collocation nodes [52], [56, 57]. It thus offers a stringent reference for accuracy and stability. On the learning-based side, we consider PINN architectures with analytical boundary and initial-condition lifting, where a neural network parameterizes a residual field multiplied by functions that enforce the constraints [40, 41]. The PDE is incorporated via a mean-squared residual loss evaluated on collocation points in space–time, and the network parameters are trained with the same optimizer and learning-rate schedule as in VSL. This allows for a balanced comparison of coefficient-space and field-space learning under closely matched numerical and optimization settings.

Within this experimental environment, VSL is applied to a set of representative PDE benchmarks that are standard in numerical analysis and machine learning for PDEs. These include elliptic problems such as the Poisson equation, parabolic problems such as the time-dependent diffusion equation, and mildly nonlinear diffusion–advection problems of Burgers type [58]. For each benchmark, we construct manufactured solutions that admit closed-form expressions for the exact solution and forcing term, enabling precise error quantification [59]. The performance metrics include relative $L^2$ and $L^\infty$ errors on test grids, as well as diagnostic quantities such as variational energy decay, PDE residual histories, initial-condition mismatch, and learning-rate trajectories. By analyzing these metrics across methods, we aim to elucidate the conditions under which VSL achieves accuracy comparable to classical solvers and identify regimes where it offers advantages over field-based learning.

VSL brings several important advantages by tightly coupling classical numerical analysis with modern machine learning practice. First, it explicitly reconnects learning-based PDE solvers with the variational viewpoint, treating the loss as a genuine space–time energy derived from strong or weak formulations rather than an ad hoc residual penalty. This alignment between the optimization objective and the underlying PDE structure enhances interpretability and provides a clear physical meaning to training progress. Second, VSL retains a high-order spectral representation of the solution, so the approximation space is mathematically transparent, compact, and capable of very rapid convergence for smooth problems, while still being trained with the same optimizers and infrastructure used for deep networks. From a numerical standpoint, spectral bases often require far fewer degrees of freedom than low-order methods or dense neural networks to reach a given accuracy, and the coefficient vector in VSL typically has moderate dimension yet represents a global, smooth approximation over the entire space–time domain. Third, because the basis functions are fixed and analytically known, the model can be regarded as a linear map from coefficients to function values, with all nonlinearity arising from the differential operator and the dependence of residuals on the solution. This turns training into a nonlinear least-squares problem with a deterministic feature map, yielding a smoother and more analyzable optimization landscape than heavily overparameterized field-space networks, which often suffer from stiffness, spectral bias, and delicate loss-balancing. Fourth, boundary and initial conditions are handled in a principled way— homogeneous Dirichlet conditions are built analytically into the basis, and initial data appear as clearly separated loss terms— so constraints are satisfied by construction rather than through fragile penalty tuning, and the approximation space is automatically contained in the correct functional setting. Finally, from the perspective of machine learning for scientific computing, VSL offers a controlled and interpretable alternative to black-box field-space neural networks: the architecture is intentionally simple—a linear model in a fixed spectral feature space trained to minimize a physics-based loss—yet, because it is implemented in standard deep learning frameworks, one can still exploit automatic differentiation, GPU acceleration, adaptive gradient methods, and advanced learning-rate schedules. In this way, VSL serves simultaneously as a strong high-order baseline, an analytically tractable laboratory for studying optimization in PDE-constrained learning, and a reusable design template for constructing new variational, coefficient-space solvers across a broad range of PDE models.

To situate the proposed framework within the broader context, the remainder of the paper is organized as follows. Section 2 reviews the relevant background and related work. Section 3 presents the spectral representation and bases, including one-dimensional Chebyshev and boundary-satisfying Dirichlet bases, and their tensor-product extension in space and time. Section 4 formalizes learning in coefficient space by defining the core VSL formulation. Section 5 discusses architectures and learning schedules, focusing on the linear coefficient model, and implementation details in modern deep learning frameworks. Sections 6, 7, and 8 report numerical experiments on three benchmark families—elliptic Poisson problems in one and two dimensions, one and two-dimensional time-dependent diffusion equation, and Burgers-type diffusion–advection problems—and compare VSL with the baseline solvers. Section 9 concludes the paper by summarizing the main findings and implications of the VSL framework for variational PDE solvers.

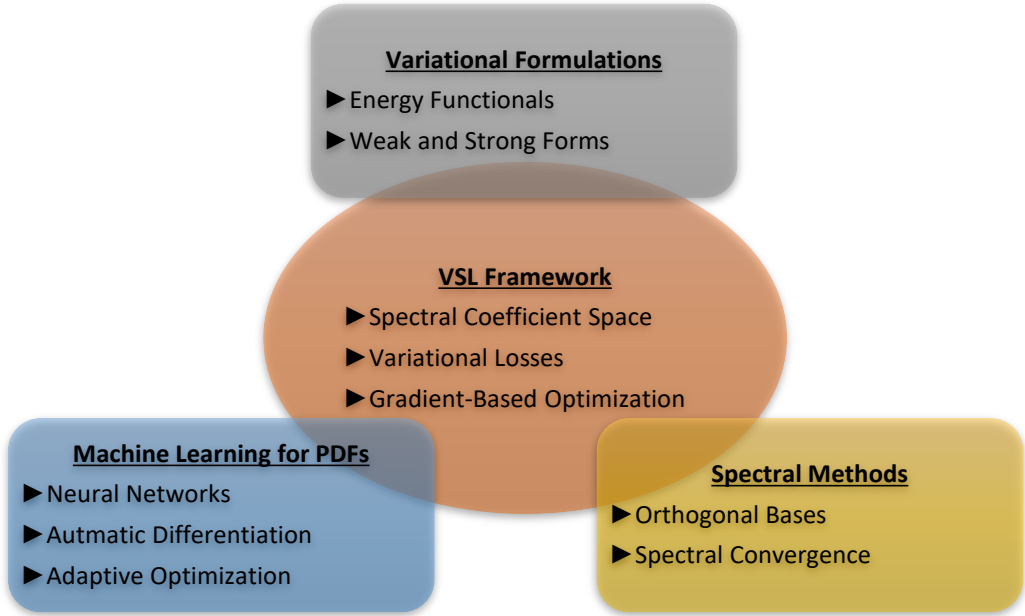

**Figure 1.** Venn diagram illustrating how variational formulations, spectral methods, and machine learning for PDEs intersect to form the VSL framework. The overlap emphasizes VSL's core idea: gradient-based optimization in spectral coefficient space using PDE-derived variational losses (strong-form least squares or Galerkin-type residuals), enabled by quadrature-based discretization and automatic differentiation within modern ML software stacks.

## 2. Background and Related Work

The study of PDEs has long been grounded in variational principles, spectral approximation theory, and, more recently, machine learning–based formulations, see figure 1. This section reviews the main strands of work that intersect in VSL: spectral methods, machine learning approaches to PDEs, and optimization strategies. The goal is not to provide an exhaustive survey, but rather to highlight the ideas and mechanisms that directly motivate and inform the VSL framework.

*2.1 Spectral methods*

Spectral methods approximate PDE solutions by global expansions in orthogonal bases, typically Fourier series on periodic domains or classical orthogonal polynomials (e.g. Chebyshev or Legendre polynomials) on bounded intervals. For a scalar field $u(x)$ on $[-1,1]$, a Chebyshev spectral approximation of degree $N$ takes the form $u_N(x) = \sum_{k=0}^{N} c_k T_k(x)$, where $T_k(x)$ are Chebyshev polynomials of the first kind and $c_k$ are expansion coefficients to be determined [11–15]. For analytic $u$, the approximation error $|u - u_N|$ typically decays at an exponential (geometric) rate in $N$, while for $u \in C^\infty$ one often observes superalgebraic decay; for finite Sobolev regularity, the decay is algebraic. This rapid decay for smooth solutions is commonly referred to as spectral (or exponential) convergence.

Several formulation variants coexist within the spectral family. In spectral Galerkin methods [11–15], [27–30], the approximate solution is inserted into the weak form, and the same global basis is used for trial and test functions, leading to systems whose entries are inner products of basis functions and their derivatives. In spectral collocation methods, the PDE is instead enforced pointwise at a set of collocation nodes, often chosen to be quadrature nodes such as Chebyshev–Gauss–Lobatto or Gauss–Legendre points. For Chebyshev collocation on $[0,1]$, one typically starts from Chebyshev–Gauss–Lobatto nodes on $[-1,1]$, $\xi_j = \cos(\pi j/N)$ for $j = 0, \ldots, N$, and maps them to $[0,1]$ via $x_j = (1 + \xi_j)/2$ [12, 13]. Differentiation matrices $D$ and $D^2$ are then constructed so that their action on nodal values approximates first and second derivatives in the physical coordinate (with the appropriate scaling induced by the affine map). The PDE then reduces to a system of ordinary differential or algebraic equations for the nodal values, subject to boundary conditions enforced by modifying rows corresponding to boundary nodes.

*2.2 Machine learning for PDEs*

Deep learning has emerged as a powerful paradigm for building surrogate models and solution operators for PDEs, with applications ranging from fluid mechanics and materials science to climate modelling and inverse problems [16, 17], [20–23]. One influential line of work is PINNs [16, 17], which parameterize the solution $u(x,t)$ by a neural network $\mathcal{N}_\theta(x,t)$ and train the parameters $\theta$ by minimizing a loss that penalizes deviations from the PDE, boundary conditions, and initial conditions at sets of collocation points [40, 41]. For a PDE of the form $F(u) = 0$, the PINN loss typically includes a term

$$\mathcal{L}_{\text{PDE}}(\theta) = \frac{1}{N_{\text{int}}} \sum_{i}^{N_{\text{int}}} |F(\mathcal{N}_\theta)(x_i, t_i)|^2, \tag{1}$$

where derivatives within $F(\mathcal{N}_\theta)$ are computed by automatic differentiation. Additional terms enforce boundary and initial conditions [16], [40, 41]. PINNs have been successfully applied to a wide range of forward and inverse PDE problems. However, training PINNs can be challenging, especially for stiff or multi-scale PDEs, high-dimensional domains, or solutions with sharp features; issues such as unbalanced loss terms, optimization pathologies, and sensitivity to initialization and sampling strategies have motivated numerous architectural and algorithmic refinements [35–39], [42, 43].

A complementary line of work targets operator learning, where the goal is to learn the mapping from PDE inputs (e.g., coefficients, forcing terms, boundary/initial data) to solution fields [20–23]. Neural operator architectures—including DeepONet-style decompositions and Fourier/neural-operator families—seek to approximate this solution operator in a data-driven manner, often enabling rapid inference once trained. These approaches are especially effective when many solves for varying inputs are required, but they typically rely on access to representative training pairs generated by numerical solvers and their generalization can depend strongly on the coverage of the training distribution. In contrast, VSL targets a physics-constrained solve for a given PDE instance by optimizing a variational objective in a fixed, interpretable approximation space, and therefore does not require supervised solution data.

A second family of approaches exploits the variational structure of PDEs more directly. In the deep Ritz method [33], for example, the PDE is expressed as the Euler–Lagrange equation of an energy functional $\mathcal{I}(u)$, and the neural network $u_\theta$ is trained by minimizing $\mathcal{I}(u_\theta)$ rather than a pointwise residual. The Deep Galerkin Method (DGM) [18] instead minimizes a strong-form residual (together with boundary/initial penalties) using stochastic sampling in space–time, thereby providing a mesh-free residual-based training procedure rather than a classical Galerkin weak formulation. Within the broader landscape of machine learning for PDEs, several hybrid approaches explicitly combine ideas from spectral methods with deep learning. For example, weak PINNs integrate the PDE residual against test functions and approximate the resulting integrals via quadrature, thereby forming losses that more closely resemble discrete Galerkin or Ritz functionals [34], [46]. Together, these developments highlight a trend towards leveraging the structure of variational PDE formulations within learning-based solvers, but most existing methods either treat the spectral or variational layers heuristically or rely on fully black-box neural architectures. VSL addresses this gap by representing the solution in a global spectral basis while defining losses that coincide with strong or weak variational energies derived from the governing PDE.

*2.3 Optimization and learning-rate schedules for PDE training*

Regardless of representation—finite elements, spectral bases, or neural networks—solving PDEs in a variational or residual-minimization framework inevitably leads to optimization problems over (possibly high-dimensional) parameter spaces [47], [51]. Classical Galerkin and spectral methods produce linear or mildly nonlinear systems that are typically solved with direct or iterative solvers such as conjugate gradients, multigrid, or Krylov subspace methods [56, 57]. Machine learning–based PDE solvers instead rely predominantly on variants of stochastic gradient descent (SGD) [60–64], in which parameters $\theta$ are updated iteratively according to $\theta_{k+1} = \theta_k - \eta_k \nabla_\theta \mathcal{L}(\theta_k)$, where $\eta_k$ is a learning rate and $\nabla_\theta \mathcal{L}(\theta_k)$ is a (possibly mini-batch) stochastic estimate of the gradient of an empirical loss [60–64]. Over the past decade, adaptive gradient methods such as Adam [48], RMSProp [61], and AdaGrad [62] have become standard defaults in deep learning due to their ability to rescale updates per-parameter based on running estimates of the first and second moments of the gradients. PINNs and related PDE-aware architectures frequently employ Adam in early training and may switch to second-order optimizers (e.g. L-BFGS) for fine-tuning, balancing the robustness of stochastic methods with the local convergence properties of quasi-Newton schemes.

Learning-rate schedules play a central role in stabilizing and accelerating optimization in these settings. Rather than keeping $\eta_k$ constant, one often employs scheduled or adaptive decay schemes in which the learning rate is decreased over the course of training or cycled between high and low values. Simple step-decay and exponential-decay schedules reduce the learning rate by a fixed factor after a prescribed number of epochs, while cyclical learning rates periodically vary $\eta_k$ between bounds to encourage exploration early on and refinement later. Cosine-annealing schedules start from a relatively large learning rate and decrease it according to a cosine curve, often combined with restarts that periodically reset the learning rate to a high value.

Taken together, the strands of work reviewed in this section motivate the design of VSL as a framework that combines the rigor and stability of classical variational formulations, the approximation power of spectral methods, and the flexibility of modern machine learning optimization. Variational formulations supply the energy functionals and weak forms that define physically meaningful losses. Spectral bases provide compact, boundary-satisfying approximation spaces with excellent convergence properties. Machine learning for PDEs contributes the tools of automatic differentiation, flexible implementation, and large-scale optimization, while the study of learning-rate schedules and adaptive methods informs the practical training of PDE-constrained models. VSL can be viewed as a coherent synthesis of these components: it performs gradient-based optimization in a spectral coefficient space using losses that are directly induced by strong-form least-squares or Galerkin-type residual constructions, with quadrature-based discretization and automatic differentiation enabling efficient and flexible implementation within modern ML software stacks.

## 3. VSL Framework: Spectral Representation and Bases

The effectiveness of the proposed VSL framework rests critically on the choice of spectral representation used to approximate the solution fields. Because the governing PDEs are posed on bounded domains in space and time and admit smooth manufactured solutions in the benchmark problems considered here, it is natural to employ global orthogonal polynomial bases that provide high-order or even spectral convergence as the number of modes increases. Among the available options, Chebyshev polynomials of the first kind offer a particularly attractive combination of numerical stability, efficient recurrence relations, and compatibility with established quadrature and differentiation techniques. In this framework, the solution is represented in a tensor-product space spanned by boundary-satisfying Chebyshev-based functions in each spatial dimension and standard Chebyshev modes in time. The resulting representation automatically enforces homogeneous Dirichlet boundary conditions in space, while maintaining sufficient flexibility in the temporal direction to encode initial conditions through an auxiliary loss. To complete the construction, space–time integrals appearing in the variational energies are evaluated via tensor-product Gauss–Legendre quadrature, which yields accurate numerical approximations without being tied to the orthogonality weight of the Chebyshev basis. This section details the underlying spectral representation, beginning with one-dimensional Chebyshev bases, then describing the construction of Dirichlet-satisfying spatial modes, extending to tensor-product bases in space and time, and finally outlining the quadrature rules used to approximate the variational energies.

In one dimension, Chebyshev polynomials of the first kind $\{T_k\}_{k \geq 0}$ are defined on the reference interval $[-1,1]$ by the relation, [11],

$$T_k(\xi) = \cos(k \arccos \xi), \qquad \xi \in [-1,1], k \in \mathbb{N}_0. \tag{2}$$

They satisfy the three-term recurrence, [12],

$$T_0(\xi) = 1, \qquad T_1(\xi) = \xi, \qquad T_{k+1}(\xi) = 2\xi T_k(\xi) - T_{k-1}(\xi), \qquad k \geq 1, \tag{3}$$

which permits stable and efficient evaluation of all modes up to a prescribed maximum degree using only multiplications and additions. Chebyshev polynomials form an orthogonal system in $L^2([-1,1], w)$ with respect to the weight $w(\xi) = (1 - \xi^2)^{-1/2}$, in the sense that, [11],

$$\int_{-1}^{1} \frac{T_j(\xi) T_k(\xi)}{\sqrt{1-\xi^2}} d\xi = \begin{cases} \pi, & j = k = 0, \\ \pi/2, & j = k \geq 1, \\ 0, & j \neq k. \end{cases} \tag{4}$$

This orthogonality, together with the ability to approximate smooth functions with near-minimal uniform error, underpins the classical theory of Chebyshev spectral method. For a smooth function $g: [-1,1] \to \mathbb{R}$, one can write a Chebyshev series, [11],

$$g(\xi) \approx \sum_{k=0}^{N} a_k T_k(\xi), \tag{5}$$

where the coefficients $a_k$ can be determined, for example, by discrete orthogonal projection using Chebyshev–Gauss or Chebyshev–Lobatto nodes [11]. For analytic functions, the approximation error typically decays at an exponential rate as $N$ increases [11–15].

To represent PDE solutions on a physical interval $[0,1]$, it is convenient to employ an affine map between the physical variable $x \in [0,1]$ and the reference coordinate $\xi \in [-1,1]$. The standard choice is

$$\xi = 2x - 1. \tag{6}$$

A function $u(x)$ defined on $[0,1]$ can thus be expanded as,

$$u(x) \approx \sum_{k=0}^{N} a_k T_k(2x - 1), \tag{7}$$

and all operations on the basis can be performed in terms of the reference variable $\xi$. In the implementation underlying the VSL framework, this mapping is used to construct Chebyshev polynomials at arbitrary points in $[0,1]$: given a vector of physical nodes $x$, one first computes $\xi = 2x - 1$ and then evaluates $T_0, \ldots, T_N$ by iterating the recurrence. This yields a feature matrix of size $M \times (N + 1)$ when $M$ evaluation points are used. Derivatives with respect to $x$ can be obtained either analytically (using derivative recurrences for $T_k$) or, as in the present work, through automatic differentiation applied to the composite mapping $(x \to \xi \to T_k(\xi))$.

For PDEs with homogeneous Dirichlet boundary conditions, it is advantageous to tailor the spatial basis so that every basis function vanishes at the domain boundaries. In one spatial dimension on $[0,1]$, the essential constraint is

$$u(0) = 0, \qquad u(1) = 0. \tag{8}$$

If a generic Chebyshev expansion is used, these conditions translate into linear constraints on the coefficients, which can be enforced by elimination or by imposing penalties. However, a more elegant and numerically robust strategy is to incorporate the boundary conditions directly into the choice of basis. To this end, one constructs a family of functions $\{\varphi_k\}_{k=0}^{N_x - 1}$ such that

$\varphi_k(0) = \varphi_k(1) = 0$ for all $k$ and $\{\varphi_k\}$ spans a rich subspace of $H_0^1(0,1)$. In the spectral representation adopted here, the Dirichlet-satisfying basis is defined as,

$$\varphi_k(x) = T_{k+2}(2x - 1) - T_k(2x - 1), \qquad k = 0, \ldots, N_x - 1. \tag{9}$$

The vanishing of $\varphi_k$ at the boundaries follows from the endpoint identities of Chebyshev polynomials. Since $T_k(1) = 1$ for all $k$ and $T_k(-1) = (-1)^k$, one has

$$\varphi_k(0) = \varphi_k(1) = 0, \tag{10}$$

so that each $\varphi_k$ satisfies the homogeneous Dirichlet condition. Moreover, the family $\{\varphi_k\}$ is linearly independent and spans a polynomial subspace that enforces zero trace at both endpoints (i.e., any linear combination satisfies $u(0) = u(1) = 0$ identically). In particular, it provides a boundary-conforming modal basis suitable for approximating functions in $H_0^1(0,1)$ with spectral accuracy for smooth targets. Any function represented as a linear combination, [11–15],

$$u(x) = \sum_{k=0}^{N_x-1} a_k \varphi_k(x), \tag{11}$$

thus automatically satisfies $u(0) = u(1) = 0$ for any coefficients $\{a_k\}$. This construction retains the smoothness and approximation power of the Chebyshev basis while "hard-coding" the essential boundary conditions into the function space, an attribute that is particularly valuable when the PDE is treated in a variational setting and the boundary values should not be controlled numerically by penalties.

The same boundary-satisfying construction is applied independently in each spatial dimension of a multi-dimensional domain. On the unit square $(0,1) \times (0,1)$, with coordinates $(x,y)$, two Dirichlet-satisfying families are introduced,

$$\varphi_i(x) = T_{i+2}(2x - 1) - T_i(2x - 1), \qquad i = 0, \ldots, N_x - 1, \tag{12.1}$$
$$\psi_j(y) = T_{j+2}(2y - 1) - T_j(2y - 1), \qquad j = 0, \ldots, N_y - 1, \tag{12.2}$$

so that $\varphi_i(0) = \varphi_i(1) = 0$ for all $i$ and $\psi_j(0) = \psi_j(1) = 0$ for all $j$. A two-dimensional spatial basis is then formed by tensor products,

$$\phi_{ij}(x,y) = \varphi_i(x)\psi_j(y), \tag{13}$$

which vanish whenever $x \in \{0,1\}$ or $y \in \{0,1\}$. The resulting span,

$$V_{N_x,N_y} = \text{span}\{\phi_{ij}(x,y)\}_{i=0,\ldots,N_x-1}^{j=0,\ldots,N_y-1} \subset H_0^1\big((0,1) \times (0,1)\big), \tag{14}$$

is a finite-dimensional subspace of the classical energy space for homogeneous Dirichlet problems. Any spectral spatial approximation constructed within $V_{N_x,N_y}$ therefore satisfies the boundary condition $u = 0$ on the entire boundary $\partial\Omega$ by design, and this property holds at all times when such a spatial expansion is coupled with a temporal basis.

In the temporal direction, the boundary conditions play a different role. For initial-value problems such as the diffusion equation, the condition $u(.,0) = u_0$ prescribes the solution at the initial time only, while no essential condition is imposed at $t = T$. Instead of enforcing the initial condition in the temporal basis, it is convenient in the present framework to use a standard Chebyshev basis in time and incorporate the initial data via an additional loss term. Let $t \in [0,T]$ denote time and introduce the affine map $\tau = 2t/T - 1 \in [-1,1]$. A temporal basis is then given by,

$$\chi_m(t) = T_m\left(\frac{2t}{T} - 1\right), \qquad m = 0, \ldots, N_t - 1. \tag{15}$$

These functions do not vanish at $t = 0$ or $t = T$; instead, they span a rich subspace of $C^0([0,T])$ capable of approximating smooth temporal behavior. The initial condition $u(.,0) = u_0(.)$ is enforced by penalizing the discrepancy between the spectral approximation and $u_0$ at $t = 0$ in the total variational energy, as described in the previous section. This separation of concerns— analytical enforcement of spatial Dirichlet constraints via boundary-satisfying bases, and numerical enforcement of the initial condition via a dedicated loss—introduces flexibility in the temporal approximation while maintaining a clean functional framework.

Combining the spatial and temporal modes yields a tensor-product space–time basis for the full solution. On the space–time cylinder $Q = (0,1)^2 \times (0,T)$, a basis element is defined as

$$\Phi_{i,j,m}(x,y,t) = \varphi_i(x)\psi_j(y)\chi_m(t), \qquad i = 0, \ldots, N_x - 1, j = 0, \ldots, N_y - 1, m = 0, \ldots, N_t - 1. \tag{16}$$

The variational spectral approximation takes the form

$$u_{N_x,N_y,N_t}(x,y,t) = \sum_{i=0}^{N_x-1}\sum_{j=0}^{N_y-1}\sum_{m=0}^{N_t-1} c_{i,j,m} \Phi_{i,j,m}(x,y,t), \tag{17}$$

where $c_{i,j,m}$ are the spectral coefficients to be learned. For convenience in implementation, the triplet index $(i,j,m)$ is flattened to a single index $n \in \{1, \ldots, N\}$, with $N = N_x N_y N_t$, and the coefficients are stored in a vector $\mathbf{c} \in \mathbb{R}^N$. Evaluating the basis at a batch of space–time points $\{(x_q, y_q, t_q)\}_{q=1}^{M_{\text{eval}}}$ yields a feature matrix $\mathbf{\Phi} \in \mathbb{R}^{M_{\text{eval}} \times N}$ whose $(q,n)$-entry is $\Phi_n(x_q, y_q, t_q)$; the approximate solution at those points is then given by the matrix–vector product

$$\mathbf{u} = \mathbf{\Phi}\mathbf{c} \in \mathbb{R}^{M_{\text{eval}}}. \tag{18}$$

In the code, this construction is implemented by first computing the one-dimensional basis matrices in $x$, $y$, and $t$, reshaping them to support broadcasting, and then forming the element-wise product to obtain the tensor-product basis. This leverages the separability of the representation and allows efficient evaluation of all basis functions at a set of points in a single operation.

The tensor-product structure is equally beneficial for the computation of spatial and temporal derivatives required in the residuals of the PDE. In the diffusion equation, for example, one needs $u_t$, $u_{xx}$, and $u_{yy}$. Analytically, derivatives of the basis functions separate as

$$\partial_t \Phi_{i,j,m}(x,y,t) = \varphi_i(x)\psi_j(y)\chi_m'(t), \qquad (19)$$
$$\partial_{xx} \Phi_{i,j,m}(x,y,t) = \varphi_i''(x)\psi_j(y)\chi_m(t), \qquad (20)$$

and similarly, for $\partial_{yy}\Phi_{i,j,m}$. Thus, once derivative bases $\varphi_i''$, $\psi_j''$ and $\chi_m'$ are available, derivatives of $u_{N_x,N_y,N_t}$ are obtained simply by replacing the corresponding factors in the tensor-product basis. Although the present implementation employs automatic differentiation to compute these derivatives numerically, the underlying tensor structure implies that the computational cost remains manageable even when higher numbers of modes are used, and it ensures smoothness and regularity of the resulting derivatives across the domain.

To evaluate the variational energies, integrals of the form

$$\int_0^T \int_0^1 \int_0^1 g(x,y,t)\, dx\, dy\, dt, \qquad (21)$$

must be approximated accurately for various integrands $g$, including squared residuals and products of residuals with basis functions. Because the basis is expressed in terms of Chebyshev polynomials, one might consider using Gauss–Chebyshev quadrature; however, the variational energies are defined with respect to the standard Lebesgue measure in space and time, not the Chebyshev weight. For this reason, the implementation adopts Gauss–Legendre quadrature on each interval and constructs a tensor-product rule for the full space–time domain. In one dimension, Gauss–Legendre quadrature on $[-1,1]$ provides nodes $\{\zeta_q\}_{q=1}^Q$ and weights $\{\omega_q\}_{q=1}^Q$ such that

$$\int_{-1}^1 h(\zeta)\, d\zeta \approx \sum_{q=1}^Q \omega_q h(\zeta_q), \qquad (22)$$

for sufficiently smooth functions $h$, [46]. Mapping to $[0,1]$ via $x = (\zeta+1)/2$ yields

$$\int_0^1 h(x)\, dx = \frac{1}{2}\int_{-1}^1 h\left(\frac{\zeta+1}{2}\right) d\zeta \approx \sum_{q=1}^Q w_q h(x_q), \qquad (23)$$

with quadrature nodes and weights given by

$$x_q = \frac{\zeta_q+1}{2}, \qquad w_q = \frac{1}{2}\omega_q. \qquad (24)$$

For the temporal interval $[0,T]$, the same Gauss–Legendre nodes $\{\zeta_q,\omega_q\}_{q=1}^{Q_t}$ on $[-1,1]$ are mapped by

$$t_q = \frac{T}{2}(\zeta_q+1), \qquad w_q^t = \frac{T}{2}\omega_q, \qquad (25)$$

so that $\int_0^T h(t)\, dt \approx \sum_{q=1}^{Q_t} w_q^t h(t_q)$.

This construction is applied independently in each dimension: Gauss–Legendre nodes and weights are generated for $x \in [0,1]$, $y \in [0,1]$, and $t \in [0,T]$. The full three-dimensional quadrature rule on $Q$ is then obtained as a tensor product of the one-dimensional rules,

$$\int_0^T \int_0^1 \int_0^1 g(x,y,t)\, dx\, dy\, dt \approx \sum_{q_x=1}^{Q_x}\sum_{q_y=1}^{Q_y}\sum_{q_t=1}^{Q_t} w_{q_x}^{(x)} w_{q_y}^{(y)} w_{q_t}^{(t)}\, g\left(x_{q_x}, y_{q_y}, t_{q_t}\right). \qquad (26)$$

In all experiments, the quadrature orders $(Q_x, Q_y, Q_t)$ are chosen sufficiently high relative to the modal degrees $(N_x, N_y, N_t)$ to avoid under-integration of squared-residual terms in the variational energies. In practice, the three-dimensional node set is constructed through a meshgrid over the one-dimensional nodes, and the combined weights are stored as the element-wise product of the one-dimensional weight arrays. This tensor-product Gauss–Legendre quadrature is used consistently for both strong-form and weak-form energies, providing high-order accuracy for smooth integrands and preserving the symmetry of the variational functionals.

In addition to quadrature nodes employed in the variational energies, the framework uses separate collocation nodes to monitor the PDE residual and to train field-space neural networks such as PINNs. For diagnostics and PINN training, it is often beneficial to sample the interior of the domain using Chebyshev–Gauss nodes, which cluster near the boundaries and thus provide enhanced resolution where boundary layers or steep gradients might arise [12]. In one dimension, Chebyshev interior nodes in $(0,1)$ are given by

$$x_j = \frac{1}{2}\left[1 - \cos\left(\frac{\pi(2j-1)}{2M_c}\right)\right], \qquad j = 1, \ldots, M_c, \tag{27}$$

which correspond to the images of Gauss–Chebyshev nodes on $[-1,1]$ under the affine map $x = (\xi + 1)/2$. These nodes are used to construct three-dimensional collocation grids in $(x, y, t)$ via tensor products, and PDE residuals evaluated at these points provide an independent measure of how well the spectral approximation satisfies the PDE in the interior. This separation between quadrature nodes for energy evaluation and collocation nodes for diagnostic residuals offers flexibility: the quadrature rules can be chosen primarily for integration accuracy, while the collocation grids can be tailored to emphasize particular regions of interest.

The combination of Chebyshev-based one-dimensional bases, Dirichlet-satisfying spatial modes, tensor-product space–time representations, and Gauss–Legendre quadrature provides a coherent and analytically transparent foundation for VSL. The Chebyshev structure ensures that smooth solutions can be represented with high accuracy using a relatively modest number of global modes, while the boundary-satisfying constructions guarantee that essential spatial boundary conditions are satisfied identically. The tensor-product construction provides a simple and efficient mechanism for extending the representation to multiple spatial dimensions and time, preserving separability and enabling compact implementations. Finally, the use of high-order quadrature rules decoupled from the orthogonality weight of the basis allows the variational energies to be defined in physically natural norms and evaluated accurately in practice. Together, these elements form the spectral backbone of the proposed framework and underpin its ability to reformulate PDE solution as a well-posed optimization problem in coefficient space.

## 4. Learning in Coefficient Space: Core VSL Formulation

Learning in coefficient space is the core mechanism by which the proposed VSL framework turns PDE solution into an optimization problem amenable to modern machine-learning tooling. Once a spectral basis has been chosen in space and time, the unknown solution field is no longer treated as a function to be discretized node-wise; instead, it is parameterized by a finite vector of spectral coefficients. The PDE, its variational energy, and any auxiliary constraints (such as initial conditions) are then expressed as differentiable functionals of these coefficients. This reformulation allows one to leverage automatic differentiation and gradient-based optimizers to compute an approximate minimizer of the variational energy directly in coefficient space. In this section, we formalize this parameterization, derive strong-form and weak-form energies as functions of the coefficient vector, incorporate initial-condition and regularization terms, and discuss the resulting gradient-based optimization problem that defines VSL.

The starting point is a finite-dimensional spectral representation of the solution. Let $Q = \Omega \times (0, T)$ denote the space–time cylinder, with $\Omega \subset \mathbb{R}^2$ a bounded spatial domain and $T > 0$ a final time. As described in the previous section, we fix a family of space–time basis functions $\{\Phi_n\}_{n=1}^N$, typically of tensor-product form $\Phi_{i,j,m}(x, y, t) = \varphi_i(x)\psi_j(y)\chi_m(t)$ in two spatial dimensions, where $\varphi_i(x)$ and $\psi_j(y)$ are Dirichlet-satisfying Chebyshev-based spatial functions and $\chi_m(t)$ are standard Chebyshev modes in time. We use a bijection $n \leftrightarrow (i, j, m)$ so $N = N_x N_y N_t$. The approximate solution is written as

$$u_N(x, t) = \sum_{n=1}^{N} c_n \Phi_n(x, t), \tag{28}$$

where, for notational simplicity, $x$ collects the spatial coordinates and $(x, t) \in Q$. The unknowns are the coefficients $\mathbf{c} = (c_1, \ldots, c_N)^T \in \mathbb{R}^N$, which form the parameter vector of the variational spectral model. The mapping

$$\mathbf{c} \mapsto u_N(\cdot, \cdot; \mathbf{c}), \tag{29}$$

is linear in $\mathbf{c}$, but the dependence of derived quantities—such as spatial derivatives, time derivatives, and nonlinear terms in the PDE—on $\mathbf{c}$ can be nonlinear because they act on $u_N$ itself. In implementation, $u_N$ is evaluated on a prescribed set of space–time points $\{z_q\}_{q=1}^M \subset Q$ (used for evaluation, diagnostics, or collocation), where $z_q = (x_q, t_q)$. Defining the basis-evaluation matrix $\mathbf{\Phi} \in \mathbb{R}^{M \times N}$ by $\Phi_{q,n} = \Phi_n(z_q)$, the vector of approximate solution values is $\mathbf{u} = \mathbf{\Phi c} \in \mathbb{R}^M$. When the points are specifically chosen as quadrature nodes for approximating integrals over $Q$, we denote their count by $M_Q$ and use the corresponding matrix $\mathbf{\Phi}_Q \in \mathbb{R}^{M_Q \times N}$ to emphasize that the discretization is quadrature-based. This representation is fully compatible with automatic differentiation frameworks, which track the dependence of $\mathbf{u}$ and its derivatives on $\mathbf{c}$ through the computational graph.

To define the variational energy, we first formalize the strong-form residual of the PDE in terms of the spectral approximation. Consider a (possibly nonlinear) PDE written abstractly as

$$\mathcal{F}(u) = 0 \text{ in } Q, \tag{30}$$

supplemented with appropriate boundary and initial conditions. In many evolution problems of interest, $\mathcal{F}$ has the structure $\mathcal{F}(u) = u_t - \mathcal{L}(u) - f$, where $\mathcal{L}$ is a spatial differential operator (e.g., diffusion or convection–diffusion) and $f$ is a forcing term [30], [55], [58]. The strong-form residual associated with the spectral approximation is

$$r(x,t;\mathbf{c}) = \mathcal{F}(u_N(\cdot,\cdot;\mathbf{c}))(x,t). \tag{31}$$

For the diffusion equation, for example, one has

$$r(x,t;\mathbf{c}) = u_{N,t}(x,t;\mathbf{c}) - \nu \Delta u_N(x,t;\mathbf{c}) - f(x,t), \tag{32}$$

where $u_{N,t}$ and $\Delta u_N$ are obtained by differentiating the spectral expansion. In practice, $r(x,t;\mathbf{c})$ is computed at a set of quadrature points $\{z_q = (x_q, t_q)\}_{q=1}^{M_Q}$ by evaluating $u_N$ and its derivatives using automatic differentiation.

The strong-form least-squares energy in coefficient space is defined by measuring the residual in an $L^2$-type norm over the space–time domain [44]. Specifically, one introduces

$$E_{\text{strong}}(\mathbf{c}) = \frac{1}{2}\int_0^T \int_\Omega [r(x,t;\mathbf{c})]^2\, dxdt, \tag{33}$$

where the factor $1/2$ is conventional and does not affect the minimizer. For smooth solutions, $E_{\text{strong}}$ is a differentiable functional of $\mathbf{c}$ because both the spectral representation and the residual are differentiable with respect to the coefficients.

In the discrete setting, all integrals over $Q$ are approximated using a tensor-product Gauss–Legendre quadrature rule. Let $\{z_q\}_{q=1}^{M_Q}$ denote the quadrature nodes on $Q$ and $\{w_q\}_{q=1}^{M_Q}$ the associated positive weights. Define the diagonal weight matrix $\mathbf{W} = \text{diag}(w_1, \ldots, w_{M_Q}) \in \mathbb{R}^{M_Q \times M_Q}$ with entries $W_{q,q} = w_q$. Collecting the basis evaluations into a feature matrix $\boldsymbol{\Phi}_Q \in \mathbb{R}^{M_Q \times N}$ with entries $(\boldsymbol{\Phi}_Q)_{q,n} = \Phi_n(z_q)$, and the residual evaluations into a vector $\mathbf{r}^Q(\mathbf{c}) \in \mathbb{R}^{M_Q}$ with components $r_q^Q(\mathbf{c}) = r(x_q, t_q; \mathbf{c})$. The discrete strong-form energy is

$$E_{\text{strong}}^Q(\mathbf{c}) \approx \frac{1}{2}\sum_{q=1}^{M_Q} w_q [r_q^Q(\mathbf{c})]^2, \tag{34.1}$$

where $M_Q = Q_x Q_y Q_t$ is the total number of quadrature points. The factor $1/2$ can be absorbed if convenient but is retained here to simplify expressions for gradients. Then the discrete strong energy is compactly

$$E_{\text{strong}}^Q(\mathbf{c}) = \frac{1}{2}[\mathbf{r}^Q(\mathbf{c})]^T \mathbf{W} \mathbf{r}^Q(\mathbf{c}). \tag{34.2}$$

The gradient of $E_{\text{strong}}^Q$ with respect to a coefficient $c_n$ is given formally by

$$\frac{\partial}{\partial c_n} E_{\text{strong}}^Q(\mathbf{c}) = \sum_{q=1}^{M_Q} w_q r(z_q; \mathbf{c}) \frac{\partial}{\partial c_n} r(z_q; \mathbf{c}), \tag{35}$$

which can be computed exactly via automatic differentiation because both $r$ and its dependence on $\mathbf{c}$ are represented in the computational graph. For linear PDEs, $\mathcal{F}$ is linear in $u$, and so $r(x,t;\mathbf{c})$ becomes an affine function of $\mathbf{c}$; in that setting, $E_{\text{strong}}$ is a quadratic functional and its minimizer can, in principle, be obtained by solving a linear system. However, using gradient-based optimization allows the same formulation to be applied uniformly to nonlinear PDEs, mixed boundary conditions, and more complex loss compositions without changing the underlying code structure.

In contrast to the strong-form energy, the weak (Galerkin) energy in coefficient space is based on projecting the residual against the spectral basis functions in an integrated sense [27–30]. For each mode index $n$, define the Galerkin residual

$$R_n(\mathbf{c}) = \int_0^T \int_\Omega r(x,t;\mathbf{c})\Phi_n(x,t)\, dxdt. \tag{36}$$

Here, "weak" refers to a Galerkin moment formulation obtained by projecting the (strong or variationally constructed) residual onto the trial space [11–15], [27–30]. When the variational form is derived by integration by parts, the resulting expressions for $R_n(\mathbf{c})$ involve lower-order derivatives (e.g., $\nabla u_N$ rather than $\Delta u_N$), improving conditioning; the present coefficient-space objective accommodates either choice through the definition of $r$ and the associated quadrature evaluation. The quantities $R_n$ measure the failure of the strong residual to be orthogonal to the finite-dimensional trial space spanned by $\{\Phi_n\}$. If $R_n(\mathbf{c}) = 0$ for all $n$, then $r(\cdot,\cdot; \mathbf{c})$ is orthogonal to the trial space, and for linear problems, this condition coincides with the classical Galerkin system for the coefficients. The weak-form variational energy is chosen as the squared $L^2$-norm of the Galerkin residual vector,

$$E_{\text{weak}}(\mathbf{c}) = \frac{1}{2}\sum_{n=1}^N [R_n(\mathbf{c})]^2 = \frac{1}{2}\|\mathbf{R}(\mathbf{c})\|_2^2, \tag{37}$$

where $\mathbf{R}(\mathbf{c}) = (R_1(\mathbf{c}), \ldots R_N(\mathbf{c}))^T \in \mathbb{R}^N$. As in the strong form, the integrals defining $R_n$ are approximated by the same Gauss–Legendre quadrature rule. Letting $\Phi_n(z_q)$ denote the $n$-th basis function evaluated at the $q$-th quadrature point, one obtains

$$R_n^Q(\mathbf{c}) \approx \sum_{q=1}^{M_Q} w_q \Phi_n(z_q) r(z_q; \mathbf{c}). \tag{38}$$

Hence, the discrete Galerkin residual can be written concisely as

$$\mathbf{R}^Q(\mathbf{c}) = \boldsymbol{\Phi}_Q^T \mathbf{W} \mathbf{r}^Q(\mathbf{c}). \tag{39}$$

The weak energy then becomes
$$E^Q_{\text{weak}}(\mathbf{c}) = \frac{1}{2}\left\|\boldsymbol{\Phi}_Q^T \mathbf{W} \mathbf{r}^Q(\mathbf{c})\right\|_2^2. \tag{40}$$
Its gradient with respect to $\mathbf{c}$ can be expressed in terms of the Jacobian of $\mathbf{r}^Q(\mathbf{c})$, but in practice it is obtained by differentiating the above expression automatically. For linear PDEs, setting $\nabla E^Q_{\text{weak}}(\mathbf{c}) = 0$ is equivalent to enforcing $R^Q_n(\mathbf{c}) = 0$ for all $n$, recovering the standard spectral Galerkin system; for nonlinear problems, the weak energy provides a natural nonlinear least-squares formulation that retains the structure of a Galerkin method while remaining compatible with gradient-based optimization.

The choice between strong-form and weak-form energies in coefficient space has important implications for numerical behavior. Strong-form energies penalize pointwise residuals equally at all quadrature nodes, which can provide good control of the PDE at specific locations but may lead to ill-conditioning when high-order derivatives are present or when the solution exhibits steep gradients. Weak-form energies, by contrast, damp residuals through integration against smooth basis functions, which can improve conditioning and more closely mimic the theoretical properties of variational formulations. In the proposed framework, both formulations are supported through a common API: the user specifies whether $E(\mathbf{c})$ is to be interpreted as $E^Q_{\text{strong}}$ or $E^Q_{\text{weak}}$, and the optimization proceeds identically in coefficient space.

For time-dependent problems, the PDE and its variational form do not automatically enforce the initial condition; this must be incorporated explicitly as an additional term in the objective. Let $u_0(x)$ denote the prescribed initial data at $t = 0$. The initial-condition discrepancy for the spectral approximation is
$$\delta_{\text{IC}}(x;\mathbf{c}) = u_N(x, 0; \mathbf{c}) - u_0(x). \tag{41}$$
An initial-condition loss is then defined as an $L^2$-type norm,
$$L_{\text{IC}}(\mathbf{c}) = \frac{1}{2}\int_\Omega |\delta_{\text{IC}}(x;\mathbf{c})|^2\, dx. \tag{42}$$
In practice, this integral is approximated using a one- or two-dimensional quadrature rule on $\Omega$, or evaluated at a set of Chebyshev interior nodes in space. In discrete form,
$$L^Q_{\text{IC}}(\mathbf{c}) = \frac{1}{2}\sum_{p=1}^{M_{\text{IC}}} \omega^{\text{IC}}_p \left[u_N(x_p, 0; \mathbf{c}) - u_0(x_p)\right]^2, \tag{43}$$
where $x_p$ are nodes used specifically for the initial time and $\omega^{\text{IC}}_p$ are associated weights. The factor $1/2$ is again included for convenience. This loss penalizes deviations of the spectral approximation from the initial data and is scaled in the total objective by a positive weight $\lambda_{\text{IC}}$ that balances its influence relative to the PDE energy.

To control the size and smoothness of the coefficient vector, and to mitigate overfitting to the particular set of quadrature and collocation points, one augments the objective with a regularization term. A simple and effective choice is an $\ell^2$-penalty on the coefficients,
$$R_{\text{reg}}(\mathbf{c}) = \frac{1}{2}\|\mathbf{c}\|_2^2 = \frac{1}{2}\sum_{n=1}^N c_n^2, \tag{44}$$
which corresponds to Tikhonov regularization in coefficient space [65]. When the underlying basis is orthogonal (or nearly so) with respect to an appropriate inner product, this penalty discourages large-amplitude, high-frequency modes and biases the solution toward smoother representations. The regularization term is scaled by a parameter $\lambda_{\text{reg}} > 0$, leading to the combined objective
$$J(\mathbf{c}) = E(\mathbf{c}) + \lambda_{\text{IC}} L_{\text{IC}}(\mathbf{c}) + \lambda_{\text{reg}} R_{\text{reg}}(\mathbf{c}), \tag{45}$$
where $E(\mathbf{c})$ is either $E_{\text{strong}}(\mathbf{c})$ or $E_{\text{weak}}(\mathbf{c})$ depending on the chosen VSL form. In discrete form, the objective used in optimization reads
$$J^Q(\mathbf{c}) = E^Q(\mathbf{c}) + \lambda_{\text{IC}} L^Q_{\text{IC}}(\mathbf{c}) + \lambda_{\text{reg}} R_{\text{reg}}(\mathbf{c}), \tag{46}$$
and is entirely expressed in terms of sums over quadrature and diagnostic nodes and algebraic functions of $\mathbf{c}$. Additional diagnostic terms, such as a PDE residual evaluated at a separate set of collocation points, can be computed to monitor solution quality but need not be included directly in the training loss.

Once the objective $J^Q(\mathbf{c})$ has been defined, learning in coefficient space reduces to an optimization problem
$$\min_{\mathbf{c}\in\mathbb{R}^N} J^Q(\mathbf{c}). \tag{47}$$
The key advantage of the spectral parameterization is that this optimization takes place in a low-dimensional Euclidean space, with $N = N_x N_y N_t$ spectral degrees of freedom typically much smaller than the number of grid points that would be required for an equivalent nodal discretization. Moreover, the mapping from $\mathbf{c}$ to the objective is fully differentiable and implemented within a deep learning framework, enabling the use of robust gradient-based optimizers designed for high-dimensional non-convex problems. In the proposed implementation, the coefficient vector $\mathbf{c}$ is declared as a trainable tensor, and the gradient

$\nabla J^Q(\mathbf{c})$ is obtained via automatic differentiation, which internally propagates sensitivities through basis evaluation, derivative computation, residual formation, and quadrature accumulation.

Gradient-based optimization proceeds iteratively. Starting from an initial guess $\mathbf{c}^{(0)}$, one generates a sequence $\{\mathbf{c}^{(k)}\}$ by, [60-64],

$$\mathbf{c}^{(k+1)} = \mathbf{c}^{(k)} - \alpha_k \mathcal{P}_k\left(\nabla J^Q(\mathbf{c}^{(k)})\right), \tag{48}$$

where $\alpha_k > 0$ is a step size and $\mathcal{P}_k$ is a preconditioning or scaling operator associated with the chosen optimizer. For simple gradient descent, $\mathcal{P}_k$ is the identity and $\alpha_k$ is a learning rate; for adaptive methods such as Adam, $\mathcal{P}_k$ rescales each component of the gradient based on estimates of its first and second moments. In the VSL framework, the Adam optimizer with a decaying or restarting learning-rate schedule is employed, providing robustness to anisotropic curvature in the objective and facilitating rapid initial progress followed by fine-tuning. Specifically, a cosine-decay-with-restarts schedule is used, in which the learning rate oscillates within progressively shrinking envelopes. This can improve early-stage exploration of the objective landscape while supporting stable fine-tuning as the effective step size decreases.

From a practical standpoint, the gradient computation is dominated by the cost of evaluating the residual and its derivatives at all quadrature points, as well as the matrix–vector multiplications required for the weak-form energy. Because the basis functions are global and smooth, relatively few quadrature points are needed to achieve high accuracy in the integrals, keeping the computational cost manageable. Furthermore, the tensor-product structure of the basis and quadrature rules allows many operations to be vectorized and executed efficiently on modern hardware accelerators. The optimization loop thus resembles that of a standard neural network training routine: at each epoch, the objective and its gradient are evaluated on a fixed set of collocation and quadrature nodes, the optimizer updates the coefficients, and diagnostics such as energy values, PDE residuals, and initial-condition errors are recorded.

The resulting learned coefficient vector $\mathbf{c}^\star$ defines a spectral approximation $u_N(\cdot,\cdot;\mathbf{c}^\star)$ that approximately minimizes the chosen variational energy subject to the initial-condition and regularization terms. Unlike neural-network-based field-space models, where the parameters are weights of a deep architecture and the solution emerges as the forward pass of a nonlinear mapping, the VSL formulation keeps the solution representation extremely compact and interpretable: each component of $\mathbf{c}^\star$ multiplies a known basis function with a clear spatial and temporal structure. This property facilitates post-processing, error analysis, and comparison with classical spectral solvers, and it highlights the central idea of VSL: to integrate the expressive power and optimization machinery of modern machine learning with the analytical clarity and convergence guarantees of spectral Galerkin methods by learning directly in the space of spectral coefficients.

## 5. Architectures and Learning Schedules

The proposed VSL framework hinges on a deliberately simple "architecture" in parameter space, combined with a nontrivial learning-rate schedule that exploits modern optimization practice. Unlike neural surrogates that parameterize the solution field through deep, nonlinear networks, VSL adopts a plain linear parameterization in the coefficients of a carefully chosen spectral basis. The complexity of the model is therefore encoded in the basis functions and in the variational energy, not in the parameterization itself. This architectural choice leads to a compact and interpretable representation that can be optimized efficiently with gradient-based methods, provided that the learning rate is controlled appropriately. In this section, the coefficient-space architecture is formalized as a linear mapping, the cosine–decay-with-restarts learning-rate schedule used in training is described in detail, and key implementation aspects are discussed, including the integration with TensorFlow's automatic differentiation, the treatment of derivative computations, and the practical configuration of optimization hyperparameters.

From an architectural viewpoint, VSL employs what may be viewed as a single-layer, purely linear network in coefficient space. Building on the coefficient-space formulation in section 4, the discrete coefficient-to-field map can be written compactly as $\mathbf{u}(\mathbf{c}) = \mathbf{\Phi}\mathbf{c}$, where $\mathbf{\Phi}$ is the fixed basis-evaluation matrix (determined by the chosen basis and node set) and $\mathbf{c}$ is the trainable coefficient vector, see figure 2. In standard neural-network notation, this is a single fully connected linear layer, $\mathbf{u} = \mathbf{W}\mathbf{c} + \mathbf{b}$, with $\mathbf{W} = \mathbf{\Phi}$ and $\mathbf{b} = \mathbf{0}$ (or a fixed lifting offset when used), i.e., no activations [66, 67] and no hidden depth [63, 64].

Expressive power is therefore not introduced through depth, width, or nonlinear activations, but through the richness of the global spectral trial space and through the physics-informed objective used for training, which may include strong-form residual norms, weak/Galerkin projections, initial- and boundary-condition penalties, and regularization. This linearity also extends to derivative computations, since spatial and temporal derivatives act on the fixed basis functions and remain linear in the coefficients, even when the governing PDE involves gradients or Laplacians. This is why the VSL model is compact and interpretable: each entry of $\mathbf{c}$ directly scales a known space–time mode, and training is the process of finding the best linear combination of these modes under the PDE-driven objective.

The learning problem remains nontrivial not because of architectural complexity, but because the loss landscape is shaped by the PDE operator and the chosen formulation; for nonlinear PDEs, composing the linear expansion with nonlinear physics terms yields residuals that are nonlinear—and typically nonconvex—in the coefficients. This framing clarifies the contrast with

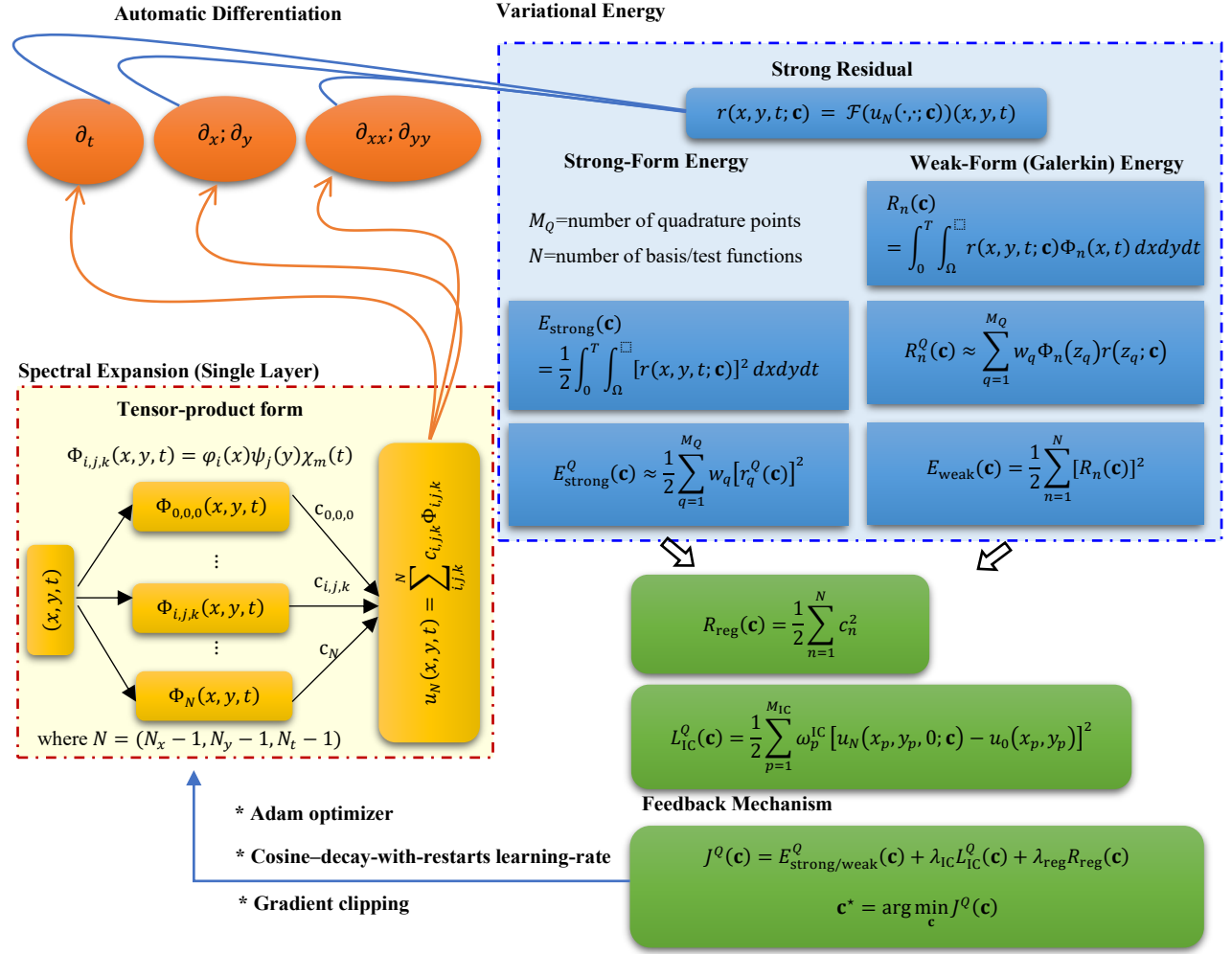

**Figure 2.** Detailed VSL architecture (coefficient-space variational solver). Flowchart of the VSL pipeline in coefficient space. A tensor-product, Dirichlet-satisfying Chebyshev basis defines a fixed feature map $\Phi_Q$ and a linear solution ansatz $u_N(z) = \Phi_Q(z)\mathbf{c}$. At Gauss–Legendre quadrature nodes $z_q$ (with weights $w_q$), automatic differentiation computes the required derivatives to assemble the strong residual $r(z_q; \mathbf{c}) = F(u_N)(z_q)$ or weak (Galerkin) moments $R^Q(\mathbf{c}) = \Phi_Q^T \mathbf{W} \mathbf{r}^Q(\mathbf{c})$. These ingredients form the discrete training objective $J^Q(\mathbf{c}) = E^Q(\mathbf{c}) + \lambda_{IC} L_{IC}^Q(\mathbf{c}) + \lambda_{reg} R_{reg}(\mathbf{c})$ (strong or weak energy, plus initial-condition and coefficient regularization terms). Optimization proceeds by gradient-based updates of the spectral coefficients $\mathbf{c}$ (e.g., Adam with cosine-decay-with-restarts and optional gradient clipping), yielding the learned coefficient vector $\mathbf{c}^\star$ and the corresponding spectral approximation $u_N(\cdot; \mathbf{c}^\star)$ over the full space–time domain.

PINNs: deep surrogates derive much of their inductive bias and expressivity from network architecture, whereas VSL encodes inductive bias explicitly through an interpretable spectral trial space and a variational/physics-based objective. Consequently, VSL provides a compact, transparent parameterization in which each coefficient scales a known mode, facilitating diagnostics and direct comparison with classical spectral methods, while enabling efficient optimization in a low-dimensional Euclidean coefficient space using standard gradient-based optimizers and modern learning-rate schedules.

For sufficiently smooth solutions, spectral expansions are known to converge rapidly, often at exponential rates with respect to the basis size. In the present framework, this translates into the ability to approximate the solution $u$ with relatively few coefficients $c_n$, especially when the basis is adapted to the boundary conditions (e.g., using Dirichlet-satisfying spatial modes) and the temporal behavior. Moreover, for linear PDEs with linear forcing, the objective function in coefficient space is quadratic, and in principle the minimizer could be obtained by solving a linear system. However, by retaining a generic gradient-based optimization loop, the same architecture can also handle nonlinear PDEs and mixed formulations without alteration. Architecturally, the VSL model thus resembles a generalized linear model with a carefully constructed, physics-informed feature space and a variational loss that encodes the PDE, boundary, and initial conditions.

The absence of hidden layers and nonlinear activation functions in the coefficient architecture also confers desirable stability properties. For instance, the mapping $\mathbf{c} \mapsto u_N$ is continuous and differentiable, and its Jacobian with respect to $\mathbf{c}$ is simply the feature matrix $\mathbf{\Phi}$. When the PDE operator is linear, the residual is an affine function of $\mathbf{c}$, and the gradient of the strong-form or weak-form energy becomes an affine transformation of $\mathbf{c}$ as well. Consequently, for linear PDE operators and quadratic loss compositions (including quadratic regularization), the coefficient-space objective is a convex quadratic function of $\mathbf{c}$, and classical results guarantee convergence of gradient descent to the unique global minimizer under appropriate step-size restrictions. In practice, we use Adam for robustness and ease of tuning, while retaining the same coefficient-space formulation for nonlinear PDEs where the objective is generally non-convex. Even when nonlinearities are present in the PDE, the structure of the basis often leads to smoother loss landscapes than those encountered in deep neural architectures, reducing the risk of pathological optimization behavior. The main challenge therefore lies not in model capacity but in efficiently navigating the coefficient space through appropriate learning-rate schedules and stopping criteria.

Given this architectural simplicity, the learning-rate schedule becomes central to practical optimization performance. In the VSL implementation, both the spectral coefficient model and the PINN baseline are trained with Adam combined with a cosine–decay-with-restarts schedule. Let $\eta_0 > 0$ denote the initial (peak) learning rate, let $\alpha \in [0,1)$ denote the minimum learning-rate fraction within each cycle, and let $k \in \mathbb{N}_0$ denote the global step index. The schedule is organized into cycles indexed by $j = 0,1,2, ...$, each of length steps, with geometrically increasing cycle lengths

$$T_j = T_0 \, (t_{\text{mul}})^j, \tag{49}$$

where $T_0 > 0$ is the first-cycle length and $t_{\text{mul}} > 0$ is a user-defined multiplier [49]. Let the cumulative step count at the start of cycle $j$ be

$$S_j = \sum_{i=0}^{j-1} T_i, \qquad S_0 = 0, \tag{50}$$

and define the normalized local step within the active cycle by

$$\tau = \frac{k - S_j}{T_j} \in [0,1), \qquad \text{for the unique } j \text{ such that } S_j \leq k < S_{j+1}. \tag{51}$$

Optionally, the peak learning rate can be scaled across cycles using a multiplier $m_{\text{mul}} > 0$, defining the cycle peak

$$\eta_j = \eta_0 (m_{\text{mul}})^j. \tag{52}$$

The cosine-decay-with-restarts learning rate at step $k$ is then

$$\eta(k) = \eta_j \left[ \alpha + \frac{1}{2}(1 - \alpha)(1 + \cos(\pi \tau)) \right], \tag{53}$$

so that $\eta(k)$ decreases smoothly from $\eta_j$ to $\alpha \eta_j$ within each cycle and restarts to $\eta_{j+1}$ at the beginning of the next cycle. In the experiments reported here, the peak rate is typically kept constant by choosing $m_{\text{mul}} = 1$.

The motivation for this schedule is twofold. First, the initial phase of each cycle, where $\tau$ is small and $\eta(k)$ is close to its peak value, enables relatively large updates that can improve early-stage exploration and reduce stagnation in flat or weakly curved regions of the objective. Second, as the cycle progresses and $\tau \to 1$, the learning rate decays towards $\eta_{\text{min}}$, effectively turning the optimizer into a fine-tuning mechanism that performs small, precise adjustments around a local basin. The restart mechanism then re-increases the learning rate at the beginning of the next cycle, allowing the optimization to explore new basins if the previous one was suboptimal. In the context of coefficient learning for PDEs, this interplay between exploration and exploitation is particularly useful when the variational energy contains multiple local minima due to nonlinearities, non-convex regularization, or complex boundary conditions.

From a numerical perspective, cosine decay with restarts also provides a principled way to reduce the learning rate over long training horizons without committing to a purely monotonic schedule. If the cycle lengths $T_j$ are chosen to grow geometrically with $j$ via $t_{\text{mul}} > 1$, early cycles are short and allow rapid exploration, while later cycles are longer and emphasize fine-tuning over extended intervals. In the VSL implementation, typical choices include $T_0$ on the order of a few hundred iterations, $t_{\text{mul}} \approx 2$, and $\alpha$ chosen such that $\eta_{\text{min}} \approx 10^{-2} \eta_0$. This configuration yields an initial phase with relatively aggressive learning rates and a terminal phase with very small steps, which is advantageous when minimizing smooth but potentially ill-conditioned objectives arising from weak-form energies and high-order derivatives.

The integration of this schedule into the VSL training loop is straightforward. The coefficient vector $\mathbf{c}$ is represented as a TensorFlow variable, and the Adam optimizer is constructed with the cosine–decay-with-restarts schedule as its learning-rate argument. At each optimization step, the current learning rate $\eta(k)$ is obtained automatically by evaluating the schedule at the optimizer's internal step counter. The training loop then proceeds by computing the loss $J^Q(\mathbf{c})$, computing its gradient with respect to $\mathbf{c}$ via automatic differentiation, optionally applying gradient clipping to limit the $\ell^2$-norm of the gradient, and finally updating $\mathbf{c}$ using the Adam update rule with the current learning rate. For diagnostic purposes, the learning rate history $\{\eta(k)\}$ is recorded in parallel with the loss values, enabling a posteriori analysis of how the schedule interacts with convergence behavior.

Gradient-based optimization in coefficient space is further stabilized by gradient clipping [68]. After computing $\nabla J^Q(\mathbf{c})$ with respect to $\mathbf{c}$, the gradient vector is clipped to have an $\ell^2$-norm not exceeding a prescribed threshold (typically of order unity). Formally, if $g = \nabla J^Q(\mathbf{c})$ and $\|g\|_2 > G_{\max}$, the clipped gradient is

$$\tilde{g} = \frac{G_{\max}}{\|g\|_2} g, \tag{54}$$

otherwise $\tilde{g} = g$. This operation prevents excessively large parameter updates that could result from transient spikes in the gradient, especially when higher-order derivatives are involved. The clipped gradient is then passed to the Adam optimizer along with the current learning rate $\eta(k)$. In TensorFlow, this operation is implemented via tf.clip_by_norm applied to the gradient with respect to the coefficient vector $\mathbf{c}$ prior to the Adam update. In the coefficient-space setting, where the number of parameters is moderate and the loss surface is relatively smooth, gradient clipping rarely activates in the later phases of training but can be useful in early epochs when the solution is far from the variational optimum.

## 6. Benchmark I – Elliptic Poisson Problems (1D and 2D)

Benchmark I focuses on elliptic Poisson problems in one and two spatial dimensions and serves as the primary testbed for evaluating the behavior of VSL on steady-state PDEs with well-understood analytical and numerical properties. The Poisson equation is chosen because it admits a canonical variational formulation, is central in potential theory and diffusion, and provides a clean setting in which to compare strong- and weak-form energies, spectral approximation, and different learning and solver paradigms. In its simplest form, the 1D benchmark considers the operator $-d^2/dx^2$ on the unit interval $\Omega = (0,1)$ with homogeneous Dirichlet boundary conditions [2]. A standard manufactured solution is prescribed as $u^\star(x) = \sin(\pi x)/\pi^2$, which leads to the right-hand side $f(x) = \sin(\pi x)$ through the relation $-u^{\star\prime\prime}(x) = f(x)$. The boundary conditions $u^\star(0) = u^\star(1) = 0$ are satisfied automatically by this choice, and the exact solution is smooth on $[0,1]$, making it ideally suited for spectral approximation. The two-dimensional extension is defined on the unit square $\Omega = (0,1)^2$ with the operator $-\Delta u = -(\partial_{xx} u + \partial_{yy} u)$ and homogeneous Dirichlet conditions on all four edges. A natural separable manufactured solution is $u^\star(x,y) = \sin(\pi x)\sin(\pi y)$, which yields $f(x,y) = 2\pi^2 \sin(\pi x)\sin(\pi y)$ and again satisfies $u^\star = 0$ on $\partial\Omega$. These 1D and 2D problems define a benchmark family in which the domain and coefficients are simple, analytical expressions for the exact solution and forcing are available, and the underlying regularity assumptions required for spectral convergence are fulfilled.

Within this benchmark, the variational structure of the Poisson equation provides a direct pathway to defining the strong and weak energies used in VSL. In the weak formulation for the 1D problem, one seeks $u \in H_0^1(0,1)$ such that

$$\int_0^1 u'(x)v'(x)dx = \int_0^1 f(x)v(x)dx, \qquad \text{for all } v \in H_0^1(0,1), \tag{55}$$

which is equivalent to the minimization of the energy functional

$$J(u) = \frac{1}{2}\int_0^1 |u'(x)|^2 dx - \int_0^1 f(x)u(x)dx. \tag{56}$$

The Euler–Lagrange equation associated with $J$ recovers the weak form, and under standard conditions the minimizer coincides with the unique solution in $H_0^1(0,1)$. An analogous structure holds in two dimensions, where the weak form reads

$$\int_\Omega \nabla u(x,y) \cdot \nabla v(x,y) dx dy = \int_\Omega f(x,y)v(x,y)dx dy, \qquad \text{for all } v \in H_0^1(\Omega), \tag{57}$$

with energy

$$J(u) = \frac{1}{2}\int_\Omega |\nabla u|^2 dx dy - \int_\Omega f\, u\, dx dy. \tag{58}$$

In addition to these classical energy functionals, one can define strong-form least-squares energies that penalize the residual of the differential operator in an $L^2$ norm. For the 1D Poisson equation, this residual is $\mathcal{R}(u)(x) = -u''(x) - f(x)$, and the strong energy is

$$E_{\text{strong}}(u) = \frac{1}{2}\int_\Omega (-u''(x) - f(x))^2 dx. \tag{59}$$

In 2D, the corresponding residual is $\mathcal{R}(u)(x,y) = -\Delta u(x,y) - f(x,y)$ and the strong energy is

$$E_{\text{strong}}(u) = \frac{1}{2}\int_\Omega (-\Delta u(x,y) - f(x,y))^2 dx dy. \tag{60}$$

In VSL, the unknown solution is not represented by a field-space neural network but rather by a spectral expansion in a boundary-satisfying basis. For the 1D problem, the approximation space is spanned by functions $\{\phi_k\}_{k=0}^{N_x-1} \subset H_0^1(0,1)$, constructed from Chebyshev polynomials mapped to the interval. This construction ensures $\phi_k(0) = \phi_k(1) = 0$ for all $k$, thereby embedding the Dirichlet boundary conditions into the approximation space. The approximate solution is written as $u_N(x; \mathbf{c}) = \sum_{k=0}^{N_x-1} c_k \phi_k(x)$, where $\mathbf{c} = (c_0, \ldots, c_{N_x-1})^T$ is the vector of spectral coefficients. In two dimensions, the approximation space is the tensor-product span of 1D bases in each coordinate. If $\{\phi_i^{(x)}\}_{i=0}^{N_x-1}$ and $\{\phi_j^{(y)}\}_{j=0}^{N_y-1}$ are boundary-satisfying bases in the $x$- and $y$-directions, the 2D basis functions are $\Phi_{ij}(x,y) = \phi_i^{(x)}(x)\phi_j^{(y)}(y)$ and the approximate solution becomes $u_N(x,y; \mathbf{c}) = \sum_{i=0}^{N_x-1}\sum_{j=0}^{N_y-1} c_{ij}\Phi_{ij}(x,y)$. The coefficient vector $\mathbf{c}$ has dimension $N_x$ in 1D and $N_x N_y$ in 2D, and the mapping from coefficients to function values is linear and analytically known.

The variational energies described above are translated into coefficient space by substituting $u_N(\cdot;\mathbf{c})$ for $u$ and approximating the relevant integrals using high-order quadrature. For the 1D strong energy (global residual), one defines

$$E_{\text{strong}}(\mathbf{c}) = \frac{1}{2}\int_0^1 \left(-u_N''(x;\mathbf{c}) - f(x)\right)^2 dx \approx \frac{1}{2}\sum_{q=1}^{Q_x} w_q \left(-u_N''(x_q;\mathbf{c}) - f(x_q)\right)^2, \tag{61}$$

where $\{x_q, w_q\}_{q=1}^{Q_x}$ are Gauss–Legendre nodes and weights on $[0,1]$, and $u_N''(x_q;\mathbf{c})$ is obtained via automatic differentiation. The 2D strong energy uses tensor-product quadrature $\{x_{q_x}, w_{q_x}\}, \{y_{q_y}, w_{q_y}\}$ on $[0,1]^2$, yielding

$$E_{\text{strong}}(\mathbf{c}) \approx \frac{1}{2}\sum_{q_x=1}^{Q_x}\sum_{q_y=1}^{Q_y} w_{q_x} w_{q_y} \left(-\Delta u_N\left(x_{q_x}, y_{q_y}; \mathbf{c}\right) - f\left(x_{q_x}, y_{q_y}\right)\right)^2. \tag{62}$$

In the weak (Galerkin) variant, one forms discrete residuals corresponding to each basis function. In 1D, the continuous Galerkin residual (the classical energy functional) for a test function $\phi_n$ is

$$R_n(\mathbf{c}) = \int_0^1 \left(u_N'(x;\mathbf{c})\phi_n'(x) - f(x)\phi_n(x)\right) dx, \tag{63}$$

which is approximated numerically via quadrature. Collecting these residuals into a vector $\mathbf{R}(\mathbf{c}) \in \mathbb{R}^{N_x}$, the weak energy is

$$E_{\text{weak}}(\mathbf{c}) = \frac{1}{2}\sum_{n=0}^{N_x-1} \left(R_n(\mathbf{c})\right)^2 = \frac{1}{2}\|\mathbf{R}(\mathbf{c})\|_2^2. \tag{64}$$

An analogous construction in 2D uses the tensor-product test functions $\Phi_{ij}$ and the standard weak-form residual

$$R_{ij}(\mathbf{c}) = \int_\Omega \left[\nabla u_N(x,y;\mathbf{c}) \cdot \nabla \Phi_{ij}(x,y) - f(x,y)\Phi_{ij}(x,y)\right] dxdy. \tag{65}$$

Leading to $E_{\text{weak}}(\mathbf{c}) = \frac{1}{2}\sum_{ij}\left(R_{ij}(\mathbf{c})\right)^2$. This corresponds to squaring the Galerkin weak-form residuals rather than projecting the strong residual directly. The resulting energy functionals $E_{\text{strong}}(\mathbf{c})$ and $E_{\text{weak}}(\mathbf{c})$ inherit the coercivity and symmetry properties of the underlying Poisson problem and provide well-structured objectives for coefficient-space learning.

The VSL setup for this benchmark consists of choosing spectral orders $N_x$ (and $N_y$ in 2D), selecting quadrature orders $Q_x$, $Q_y$ that are sufficient to integrate the relevant polynomials accurately, and initializing the coefficient vector $\mathbf{c}$, typically at zero. The coefficients are treated as trainable parameters and are updated by minimizing $E_{\text{strong}}(\mathbf{c})$ or $E_{\text{weak}}(\mathbf{c})$ via gradient-based methods. Gradients $\nabla_{\mathbf{c}} E_{\text{strong}}$ and $\nabla_{\mathbf{c}} E_{\text{weak}}$ are computed using reverse-mode automatic differentiation, which handles the nested differentiation required to obtain second derivatives and gradient contributions efficiently. For the elliptic Poisson benchmarks, there is no time dimension and no initial condition term, so the total loss in VSL reduces to the chosen variational energy plus an optional $\ell^2$-regularization term $\lambda_{\text{reg}}\|\mathbf{c}\|_2^2$ that can stabilize optimization and filter high-frequency modes. In practice, the loss minimized in VSL takes the form

$$J_{\text{VSL}}(\mathbf{c}) = E_{\text{strong/weak}}(\mathbf{c}) + \lambda_{\text{reg}}\|\mathbf{c}\|_2^2, \tag{66}$$

with the choice between strong and weak formulation explored as part of the benchmark.

To contextualize the performance of VSL, two families of baselines are considered: classical spectral collocation solvers and field-space PINNs with analytic boundary lifting. The classical 1D collocation method employs Chebyshev–Lobatto nodes $\{x_j\}_{j=0}^{N_{\text{spec}}}$ on $[0,1]$, constructed by mapping nodes $z_j = \cos(\pi j/N_{\text{spec}})$ from $[-1,1]$ to the physical interval. A differentiation matrix $D$ approximating $d/dx$ is constructed using standard Chebyshev formulas, and the second-derivative matrix $D^2$ approximates $d^2/dx^2$. By restricting $D^2$ to interior nodes and enforcing $u = 0$ on the boundary, one obtains a linear system $-D_{\text{int}}^2 u_{\text{int}} = f_{\text{int}}$, where $u_{\text{int}}$ and $f_{\text{int}}$ are the interior values of the solution and forcing. This system is solved by direct factorization. In two dimensions, Chebyshev–Lobatto grids in $x$ and $y$ are combined, and a discrete Laplacian is assembled via a Kronecker sum of 1D second-derivative matrices, again with boundary conditions enforced by eliminating boundary nodes. The resulting algebraic system $-Lu_{\text{int}} = f_{\text{int}}$ provides a high-order reference solution from which errors can be measured.

The PINN baselines approximate the same Poisson problems using shallow feedforward neural networks in field space with analytic lifting to enforce boundary conditions. In 1D, the solution is parameterized as

$$u_\theta(x) = x(1-x)v_\theta(x), \tag{67}$$

where $v_\theta$ is a neural network with parameters $\theta$ and smooth non-linear activations. The multiplicative factor $x(1-x)$ ensures $u_\theta(0) = u_\theta(1) = 0$, so no boundary penalty is required. The loss is defined as a mean-squared PDE residual evaluated at collocation points,

$$\mathcal{L}_{\text{PINN}}(\theta) = \frac{1}{|\mathcal{S}|} \sum_{x_i \in \mathcal{S}} \left(-u_\theta''(x_i) - f(x_i)\right)^2, \tag{68}$$

with derivatives obtained via automatic differentiation. In 2D, a similar lifting is used,

$$u_\theta(x, y) = x(1-x)y(1-y)v_\theta(x, y), \tag{69}$$

and the loss penalizes the mean-squared residual of $-\Delta u_\theta - f$ on a tensor-product collocation set. The PINNs are trained with the same optimizer and learning-rate schedule as VSL to isolate the effect of parameterization and loss definition rather than differences in optimization settings.

The results of the Poisson benchmarks are assessed using both discretization-agnostic and method-specific metrics. For each method, the approximate solution $u_{\text{num}}$ is evaluated on a dense uniform test grid and compared to the exact solution $u^\star$. The relative $L^2$ error, [59],

$$L^2_{\text{rel}} = \|u_{\text{num}} - u^\star\|_{L^2} / \|u^\star\|_{L^2}, \tag{70}$$

and the relative $L^\infty$ error,

$$L^\infty_{\text{rel}} = \|u_{\text{num}} - u^\star\|_{L^\infty} / \|u^\star\|_{L^\infty}, \tag{71}$$

are computed numerically via quadrature or discrete norms on the test grid. For VSL, one additionally monitors the decay of the variational energy $E_{\text{strong/weak}}(\mathbf{c})$ over training epochs, as well as any diagnostic residual evaluated at independent collocation points. For the classical collocation solver, the focus is on the discretization error as a function of spectral resolution $N_{\text{spec}}$. For PINNs, training trajectories of the residual loss and potential overfitting to collocation points versus generalization to independent test grids are examined.

In this benchmark for the 1D Poisson boundary-value problem, the three solvers exhibit clearly distinct accuracy–cost profiles (see figure 3). The classical Chebyshev–Lobatto collocation method, with $N_{\text{spec}} = 32$, is the most accurate and efficient: on the dense test grid it achieves $L^2_{\text{rel}} = 8.429 \times 10^{-8}$ and $L^\infty_{\text{rel}} = 2.114 \times 10^{-7}$. VSL in coefficient space performs competitively when trained with the weak-form energy functional. Using a Dirichlet Chebyshev basis with $N = 16$ nodes only, the energy-based run reaches near-spectral accuracy with $L^2_{\text{rel}} = 5.025 \times 10^{-7}$ and $L^\infty_{\text{rel}} = 5.883 \times 10^{-7}$. Moreover, the diagnostic strong-form residual (evaluated on interior Chebyshev points) decreases from $\sim 5.23 \times 10^2$ to the prescribed tolerance $10^{-10}$, stopping at epoch $\approx 941$, indicating that the learned coefficient vector corresponds to a high-quality weak-form minimizer that also satisfies the PDE strongly to high precision. The hard-boundary-condition PINN (depth 4 and width 64 with tanh activations) converges robustly without boundary-weight tuning (the boundary loss remains identically zero by construction) and steadily reduces the PDE residual loss to $\sim 2 \times 10^{-5} - 10^{-4}$ over 3000 epochs, yielding a reasonable but less accurate approximation with $L^2_{\text{rel}} = 6.218 \times 10^{-4}$ and $L^\infty_{\text{rel}} = 8.791 \times 10^{-4}$. Overall, for this smooth linear problem the spectral solvers (collocation and energy-based VSL) deliver near-machine-precision solutions, with collocation being dramatically more efficient, while the PINN provides a stable physics-constrained approximation but remains less competitive in accuracy for a comparable training budget.

When the VSL objective is switched from the variational energy (weak form) to the integrated strong-form residual $E_{\text{strong}}(\mathbf{c}) \approx \int_0^1 (u_{xx} + \sin(\pi x))^2 dx$, the training behavior changes markedly and the final accuracy deteriorates (see figure 4). The global-residual VSL solution attains only $L^2_{\text{rel}} = 4.420 \times 10^{-4}$ and $L^\infty_{\text{rel}} = 7.374 \times 10^{-4}$, i.e., comparable to (or slightly better than) the PINN but far inferior to the energy-based VSL. This contrast highlights that, for this Poisson BVP, the weak-form energy functional provides a significantly better-conditioned and more stable optimization objective in coefficient space, whereas directly minimizing the strong-form residual—requiring second derivatives and thus being more sensitive to float32 automatic differentiation and aggressive warm restarts—leads to oscillatory training and inferior accuracy for a substantially larger computational budget.

Note that, because both VSL and the PINN are trained with stochastic gradient-based optimization under float32 automatic differentiation, the reported stopping epochs, loss trajectories, and final error metrics should be interpreted as representative outcomes from a single run rather than deterministic constants; in repeated executions, modest run-to-run variability is expected in the final $L^2_{\text{rel}}$ and $L^\infty_{\text{rel}}$ errors as well as in the epoch at which the diagnostic tolerance is reached (particularly for the VSL early-stopping criterion and for the non-monotone strong-form/global-residual objective).

In the 2D Poisson benchmark with homogeneous Dirichlet boundary conditions and exact solution $u^\star(x, y) = \sin(\pi x)\sin(\pi y)$, the three solvers exhibit clearly different accuracy–cost trade-offs on the common Chebyshev–Lobatto evaluation grid $N_{\text{spec}} = 24$, (see figures 5 and 6). The classical 2D Chebyshev collocation method with $(N_x = N_y = 25)$ is essentially exact for this smooth linear problem: it reaches near machine precision with relative errors $L^2_{\text{rel}} \approx 8.09 \times 10^{-15}$ and $L^\infty_{\text{rel}} \approx 9.44 \times 10^{-15}$.

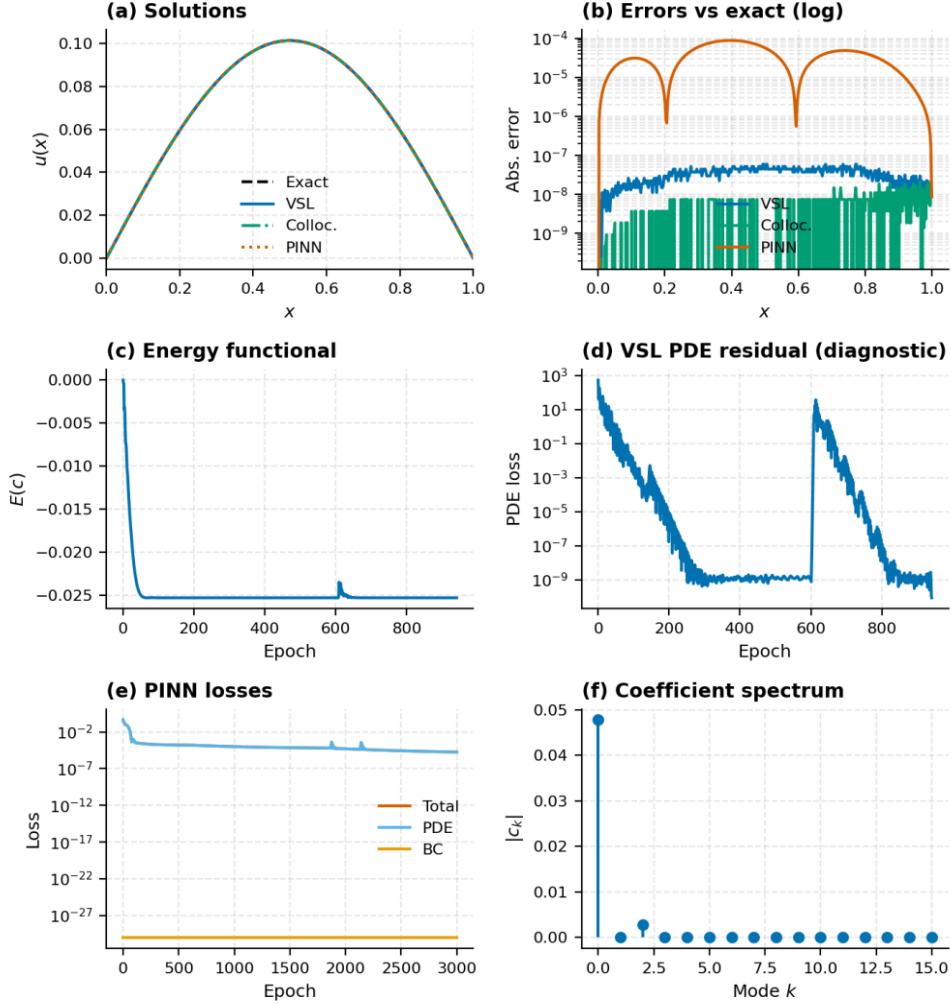

**Figure 3.** Overview comparison (3×2 grid) for the 1D Poisson boundary-value problem $u''(x) = -\sin(\pi x)$ on $x \in [0,1]$ with homogeneous Dirichlet conditions $u(0) = u(1) = 0$. Panel (a) compares the dense-grid solutions for the exact solution, VSL trained by the weak-form Poisson energy (Dirichlet Chebyshev basis, $N = 16$), classical Chebyshev–Lobatto collocation ($N_{\text{spec}} = 32$), and a hard-BC PINN (tanh MLP, width 64, depth 4). Panel (b) shows the pointwise absolute error $|u_{\text{num}} - u_{\text{exact}}|$ on a logarithmic scale. Panels (c)–(d) report the VSL training history: (c) the primary objective (energy functional) versus epoch and (d) the diagnostic strong-form residual loss $\langle (u_{xx} + \sin(\pi x))^2 \rangle$ evaluated at Chebyshev interior points. Panel (e) shows the PINN training losses (PDE loss, with the boundary loss shown only for diagnostics because the hard-BC ansatz enforces $u(0) = u(1) = 0$ exactly). Panel (f) displays the learned VSL coefficient spectrum $|c_k|$ versus mode index $k$.

When VSL is trained in coefficient space using the standard weak-form energy functional, the method attains near-spectral accuracy within a compact tensor-product Dirichlet Chebyshev basis ($N_x = N_y = 8$). In the reported run, the VSL solution achieves $L^2_{\text{rel}} \approx 1.38 \times 10^{-7}$ and $L^\infty_{\text{rel}} \approx 1.79 \times 10^{-7}$. The training curves show rapid convergence of both the quadrature-weighted strong-form residual and the diagnostic collocation residual to $\mathcal{O}(10^{-8})$ by roughly $\sim 400$ epochs, consistent with the fact that minimizing $E_{\text{weak}}(u)$ provides a well-conditioned objective for Poisson problems and avoids directly optimizing second-derivative residuals. The boundary-lifted PINN baseline (depth 4 and width 64 with tanh activations) converges more slowly in accuracy but remains robust without any boundary-weight tuning. In the weak-form campaign, the PINN reduces its PDE residual from $\mathcal{O}(10^1)$ to $\mathcal{O}(10^{-4})$ over the training window and yields a moderate-accuracy solution with $L^2_{\text{rel}} \approx 9.27 \times 10^{-5}$ and $L^\infty_{\text{rel}} \approx 1.39 \times 10^{-4}$, i.e. less accurate than energy-based VSL.

Switching VSL from the variational energy objective to the global strong-form least-squares objective changes the training behavior markedly (see figures 7 and 8). Although the quadrature residual decreases steadily early in training (e.g. from $\mathcal{O}(10^2)$ toward $\mathcal{O}(10^{-2})$), the optimization becomes significantly less stable under cosine warm restarts: the objective and the diagnostic collocation residual show pronounced non-monotone behavior and repeated increases after restart events (for example, near epoch $\sim 500$ and again near $\sim 1500$, where the global residual jumps back to $\mathcal{O}(10^{-1}) - \mathcal{O}(10^1)$).

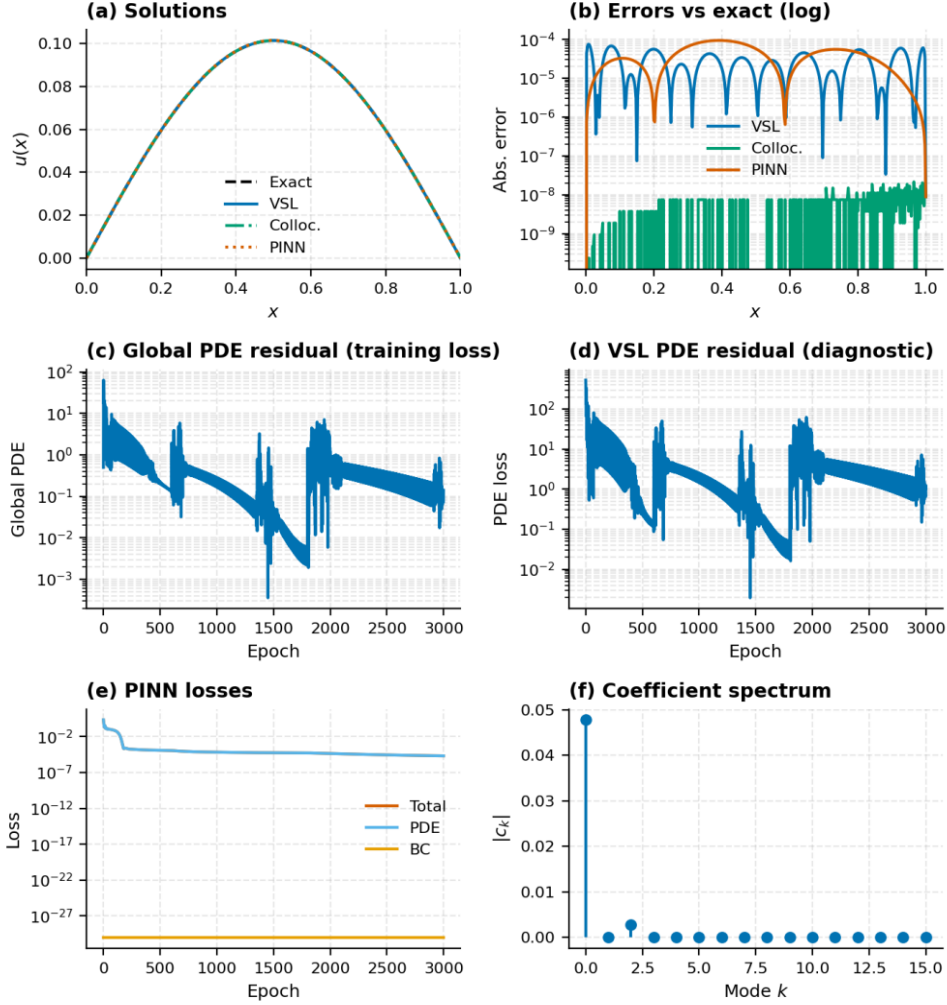

**Figure 4.** Overview comparison (3×2 grid) for the same Poisson benchmark, with VSL trained by the integrated strong-form residual. The panels match figure 3: (a) dense-grid solution comparison (exact vs global-PDE VSL with Dirichlet Chebyshev basis $N = 16$ vs collocation $N_{\text{spec}} = 32$ vs hard-BC PINN), (b) pointwise absolute errors (log scale), (c) the VSL primary objective $E_{\text{strong}}(\mathbf{c}) \approx \int_0^1 (u_{xx} + \sin(\pi x))^2 dx$ computed by Gauss–Legendre quadrature, (d) the diagnostic residual loss evaluated at Chebyshev interior points, (e) PINN loss histories, and (f) the learned VSL coefficient spectrum $|c_k|$. This figure highlights the increased non-monotonicity and reduced accuracy when replacing the weak-form energy objective by the global strong-form residual objective.

As a consequence, the final accuracy degrades substantially relative to the weak-form run: the global strong-form VSL solution ends at $L^2_{\text{rel}} \approx 8.33 \times 10^{-5}$ and $L^\infty_{\text{rel}} \approx 1.35 \times 10^{-4}$. For comparison, in the same strong-form campaign the PINN achieves $L^2_{\text{rel}} \approx 6.64 \times 10^{-5}$ and $L^\infty_{\text{rel}} \approx 7.12 \times 10^{-5}$, making it competitive with (and slightly better than) the global strong-form VSL result.

The analysis of these benchmarks emphasizes several qualitative aspects that are directly relevant to the design of VSL. First, because the manufactured solutions are smooth and the domain is simple, the approximation error associated with the spectral bases in VSL and classical collocation decays rapidly with increasing basis size, consistent with spectral convergence theory. Once optimization has reduced the variational energy to values near machine precision, the dominant error in VSL is expected to be this spectral approximation error, and thus VSL inherits the convergence behavior of the underlying spectral space. Second, the strong- and weak-form energies provide complementary views of the solution quality. In the strong form, minimizing $\int (-\Delta u_N - f)^2$ enforces pointwise balance of the PDE up to quadrature errors, while in the weak form, minimizing $\sum (R_n(\mathbf{c}))^2$ enforces orthogonality of the residual against the chosen basis. For the Poisson operator, both energies are well conditioned, and their minimizers coincide with the Galerkin solution under exact quadrature, so differences observed in practice can be attributed to quadrature and optimization details. Third, comparing VSL to field-space PINNs highlights the impact of parameterization and loss structure. VSL operates in a low-dimensional coefficient space with a linear mapping to function values and a loss that is strictly derived from variational principles, whereas PINNs operate in a much higher-dimensional weight space with a more complex loss landscape. This difference manifests in the structure of optimization trajectories: VSL's energy decay is typically monotone and aligned with physically meaningful quantities, while PINN loss curves can exhibit more irregular behavior.

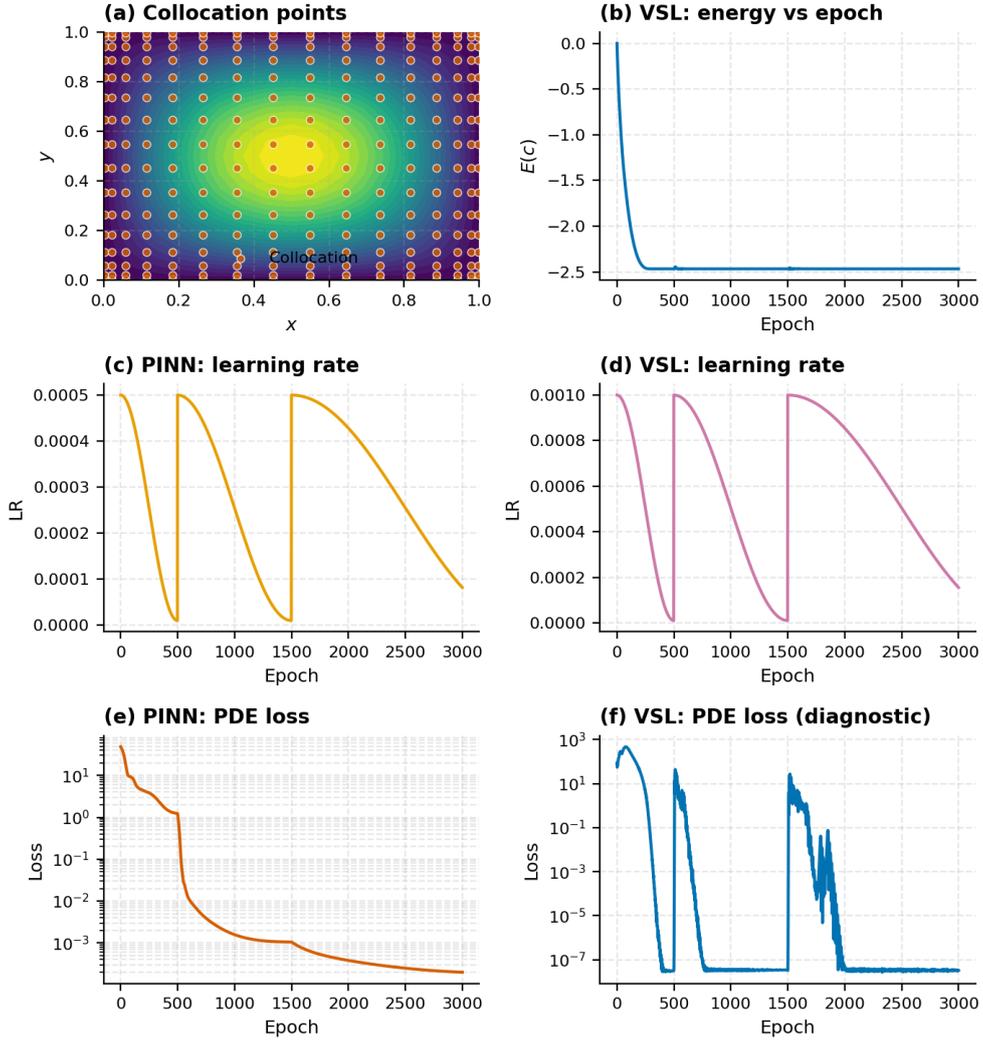

**Figure 5.** Training diagnostics for the 2D Poisson benchmark weak form objective. (a) Interior collocation points used for the diagnostic PDE residual and for PINN training, shown over the exact solution contours. (b) VSL training objective versus epoch, where the optimized objective is the discrete energy $E(\mathbf{c})$ evaluated by tensor-product Gauss–Legendre quadrature. (c) PINN learning-rate schedule (cosine annealing with warm-restart parameters. (d) VSL learning-rate schedule (cosine annealing with the same restart parameters and base learning rate $10^{-3}$). (e) PINN PDE residual loss versus epoch, computed as a mean-squared strong-form residual $(u_{xx} + u_{yy} + f)$ at collocation points under boundary lifting. (f) VSL diagnostic PDE residual loss versus epoch, computed on the collocation grid as the mean-squared strong-form residual; in the energy-based setting this diagnostic is not the optimized objective but serves to verify that the learned coefficient vector yields a solution consistent with the PDE.

## 7. Benchmark II – 1D and 2D Time-Dependent Diffusion Equation

Benchmark II considers one- and two-dimensional time-dependent diffusion equation and is designed to assess the behavior of VSL on a prototypical parabolic problem in which both spatial and temporal dynamics must be captured with high fidelity. The model problem is posed on the space–time cylinder $Q = \Omega \times I$, where $\Omega = (0,1)^2$ is the unit square in $\mathbb{R}^2$ and $I = [0,1]$ is the unit time interval. The governing equation is the heat-type diffusion equation [30]

$$u_t(x,y,t) - \nu\left(u_{xx}(x,y,t) + u_{yy}(x,y,t)\right) = f(x,y,t), \quad (x,y,t) \in Q, \tag{72}$$

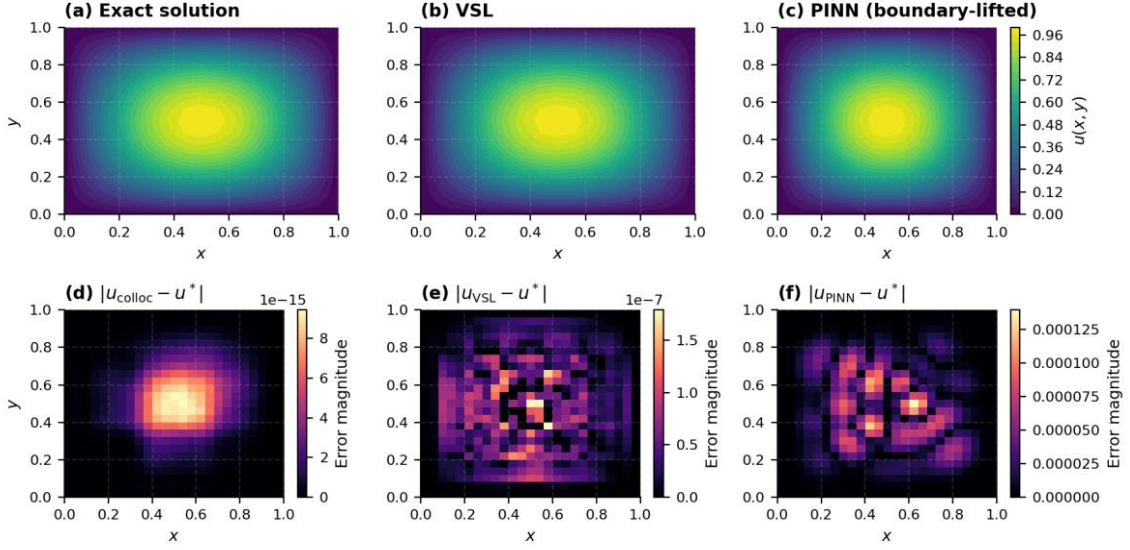

**Figure 6.** 2D Poisson: solutions and error fields (weak form objective). Top row: (a) exact solution $u^\star(x,y) = \sin(\pi x)\sin(\pi y)$, (b) VSL solution $u_{\text{VSL}}$ represented in a tensor-product Dirichlet Chebyshev basis with $(N_x, N_y) = (8,8)$, and (c) boundary-lifted PINN solution $u_{\text{PINN}} = x(1-x)y(1-y)v_\theta(x,y)$. Bottom row: absolute error fields on the common Chebyshev–Lobatto evaluation grid: (d) $|u_{\text{colloc}} - u^\star|$ from classical 2D Chebyshev collocation (Dirichlet enforced on the boundary and linear interior solve), (e) $|u_{\text{VSL}} - u^\star|$, and (f) $|u_{\text{PINN}} - u^\star|$.

with a constant diffusion coefficient $\nu > 0$. Homogeneous Dirichlet boundary conditions are prescribed along the spatial boundary,

$$u(x,y,t) = 0, \qquad (x,y) \in \partial\Omega, t \in I, \tag{73}$$

and the initial condition is given by

$$u(x,y,0) = u_0(x,y) = \sin(\pi x)\sin(\pi y), \qquad (x,y) \in \Omega. \tag{74}$$

To enable precise error quantification, a manufactured exact solution is adopted in the form

$$u^\star(x,y,t) = e^{-t}\sin(\pi x)\sin(\pi y), \tag{75}$$

which is smooth in both space and time and compatible with the boundary conditions. Substituting $u^\star$ into the PDE yields

$$u_t^\star = -e^{-t}\sin(\pi x)\sin(\pi y), \qquad u_{xx}^\star + u_{yy}^\star = -2\pi^2 e^{-t}\sin(\pi x)\sin(\pi y), \tag{76}$$

so the forcing term is

$$f(x,y,t) = (2\nu\pi^2 - 1)e^{-t}\sin(\pi x)\sin(\pi y). \tag{77}$$

This construction ensures that $u^\star$ satisfies the PDE, the homogeneous Dirichlet boundary conditions, and the initial condition exactly, providing an analytically tractable benchmark for comparing numerical solvers.

For the corresponding 1D diffusion benchmark on $\Omega = (0,1)$ and $t \in (0,1)$, we use

$$u_t(x,t) - \nu u_{xx}(x,t) = f(x,t), \qquad u(0,t) = u(1,t) = 0, \qquad u(x,0) = \sin(\pi x), \tag{78}$$

with manufactured solution $u^\star(x,t) = e^{-t}\sin(\pi x)$. Substitution yields $u_t^\star = -e^{-t}\sin(\pi x)$ and $u_{xx}^\star = -\pi^2 e^{-t}\sin(\pi x)$, hence

$$f(x,t) = (\nu\pi^2 - 1)e^{-t}\sin(\pi x). \tag{79}$$

We take the 2D diffusion benchmark as the primary test problem and present its results first. The 1D diffusion benchmark is then reported as a simpler special case, obtained by a spatial reduction of the 2D setting.

The VSL implementation introduces a finite-dimensional spectral ansatz in coefficient space for $u(x,y,t)$, constructed from Dirichlet-satisfying Chebyshev bases in space and a Chebyshev basis in time. In each spatial direction, say $x \in [0,1]$, a family of basis functions $\{\phi_k^{(x)}\}_{k=0}^{N_x-1} \subset H_0^1(0,1)$ is defined via a mapped Chebyshev construction of the form $\phi_k^{(x)}(x) = T_{k+2}(2x - 1) - T_k(2x - 1)$, where $T_n$ is the $n$-th Chebyshev polynomial of the first kind. This definition ensures $\phi_k^{(x)}(0) = \phi_k^{(x)}(1) = 0$ for all $k$, so homogeneous Dirichlet boundary conditions in the $x$-direction are built into the approximation space. An analogous family $\{\phi_\ell^{(y)}\}_{\ell=0}^{N_y-1}$ is constructed in the $y$-direction. In time, one employs a Chebyshev basis $\{\psi_m\}_{m=0}^{N_t-1}$ on $[0,1]$ defined by $\psi_m(t) = T_m(2t - 1)$. The full space–time approximation space is then the tensor-product span of these 1D bases, with basis functions $\Phi_{k\ell m}(x,y,t) = \phi_k^{(x)}(x)\phi_\ell^{(y)}(y)\psi_m(t)$, for indices $k = 0, \ldots, N_x - 1$, $\ell = 0, \ldots, N_y - 1$, and $m = 0, \ldots, N_t - 1$. The approximate solution is written as

$$u_N(x, y, t; \mathbf{c}) = \sum_{k=0}^{N_x-1} \sum_{\ell=0}^{N_y-1} \sum_{m=0}^{N_t-1} c_{k\ell m} \Phi_{k\ell m}(x, y, t), \tag{80}$$

where $\mathbf{c}$ denotes the vector of spectral coefficients of dimension $N_x N_y N_t$. Because each spatial basis function vanishes on the boundary, any coefficient vector $\mathbf{c}$ yields an approximation $u_N$ that satisfies the homogeneous Dirichlet conditions exactly, without the need for penalty terms or constraint enforcement during optimization.

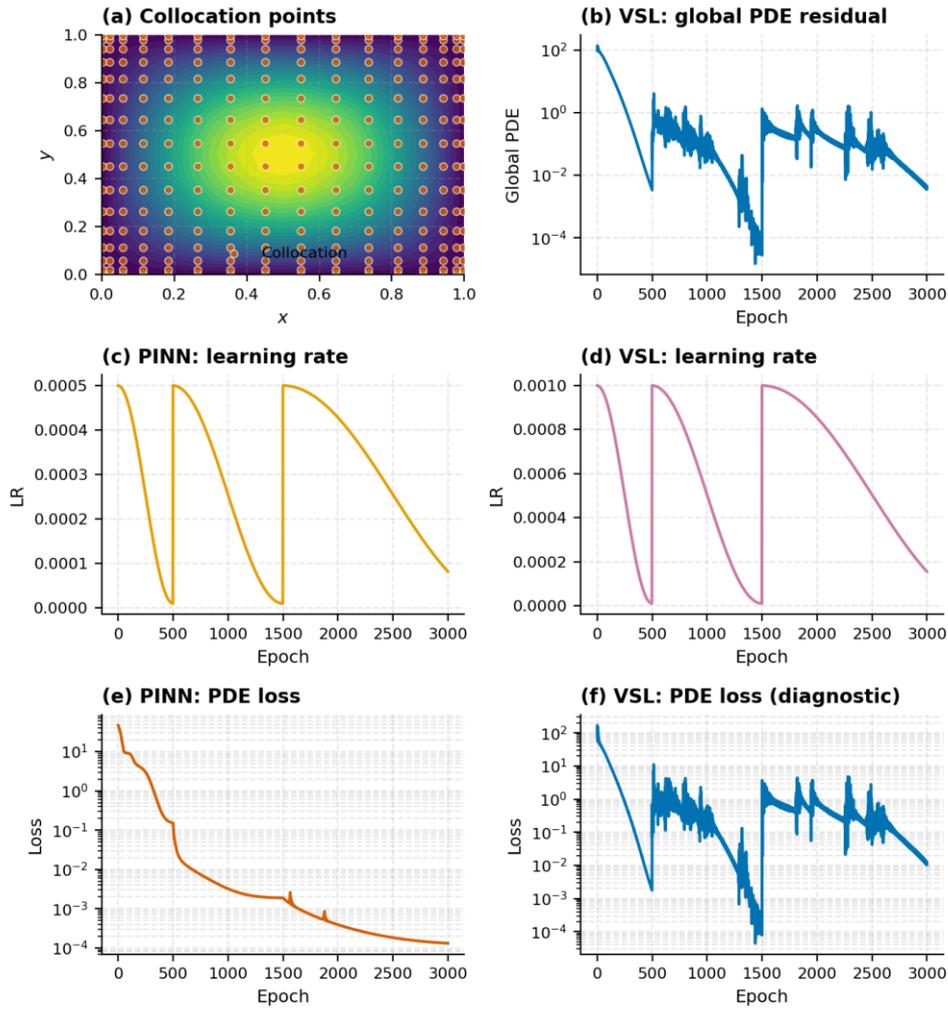

**Figure 7.** 2D Poisson: The same as in figure 5, but for VSL strong-form objective.

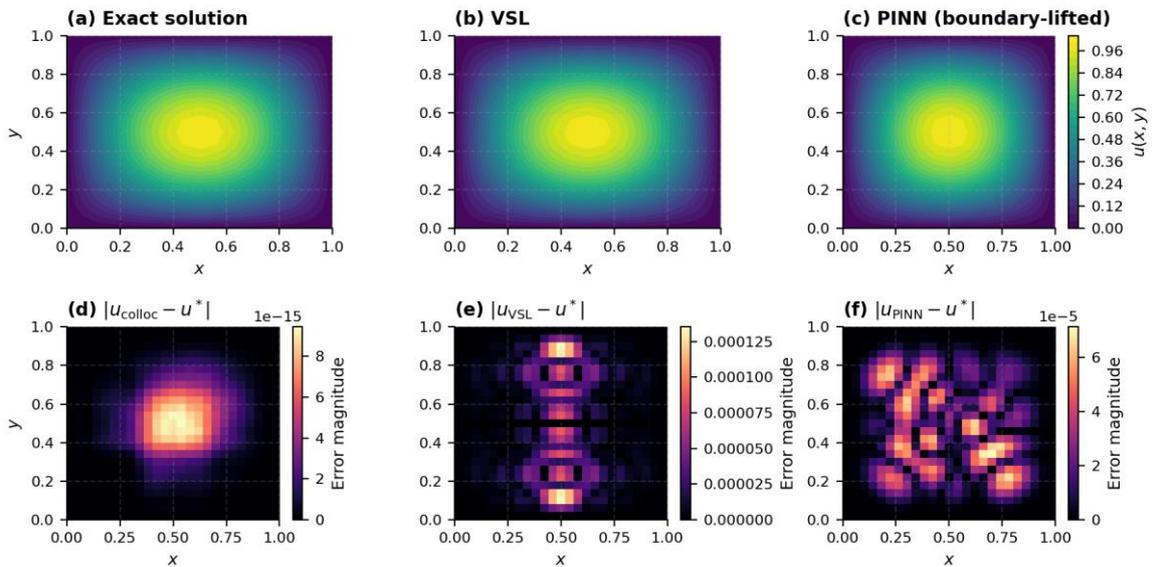

**Figure 8.** 2D Poisson: The same as in figure 6, but for VSL strong-form objective.

The strong-form VSL energy in coefficient space for the diffusion problem is obtained by substituting $u_N(\cdot; \mathbf{c})$ into $E_{\text{strong}}(u)$ and approximating the space–time integral by tensor-product Gauss–Legendre quadrature. Let $\{x_{q_x}, w_{q_x}\}$, $\{y_{q_y}, w_{q_y}\}$, and $\{t_{q_t}, w_{q_t}\}$ be Gauss–Legendre nodes and weights on $[0,1]$ in each coordinate, with $q_x = 1, \ldots, Q_x$, $q_y = 1, \ldots, Q_y$, and $q_t = 1, \ldots, Q_t$. The quadrature approximation of the strong energy is

$$E_{\text{strong}}(\mathbf{c}) \approx \frac{1}{2} \sum_{q_x=1}^{Q_x} \sum_{q_y=1}^{Q_y} \sum_{q_t=1}^{Q_t} w_{q_x} w_{q_y} w_{q_t} \left( \mathfrak{M}\left( x_{q_x}, y_{q_y}, t_{q_t}; \mathbf{c} \right) \right)^2, \tag{81}$$

$$\mathfrak{M}\left( x_{q_x}, y_{q_y}, t_{q_t}; \mathbf{c} \right) = u_{N,t}\left( x_{q_x}, y_{q_y}, t_{q_t}; \mathbf{c} \right) - \nu \left( u_{N,xx}\left( x_{q_x}, y_{q_y}, t_{q_t}; \mathbf{c} \right) + u_{N,yy}\left( x_{q_x}, y_{q_y}, t_{q_t}; \mathbf{c} \right) \right) - f\left( x_{q_x}, y_{q_y}, t_{q_t} \right). \tag{82}$$

Spatial and temporal derivatives of $u_N$ are computed using nested automatic differentiation tapes, eliminating the need to derive analytical derivative expressions for each basis function. The weak-form energy, moment-residual energy or Galerkin least-squares (GLS), is constructed by first forming Galerkin residuals

$$R_{k\ell m}(\mathbf{c}) = \iiint_Q \mathcal{R}(u_N(\cdot; \mathbf{c}))(x, y, t) \Phi_{k\ell m}(x, y, t) \, dx \, dy \, dt, \tag{83}$$

where $\mathcal{R}(u_N(\cdot; \mathbf{c}))(x, y, t) = u_{N,t}(x, y, t; \mathbf{c}) - \nu \left( u_{N,xx}(x, y, t; \mathbf{c}) + u_{N,yy}(x, y, t; \mathbf{c}) \right) - f(x, y, t)$, approximated via the same tensor-product quadrature, and then defining

$$E_{\text{weak(GLS)}}(\mathbf{c}) = \frac{1}{2} \sum_{k=0}^{N_x-1} \sum_{\ell=0}^{N_y-1} \sum_{m=0}^{N_t-1} \left( R_{k\ell m}(\mathbf{c}) \right)^2. \tag{84}$$

Importantly, no integration by parts in space is performed inside the energy. That is, the differential operator remains acting on the trial field $u_N$ inside $\mathcal{R}(u_N)$, and the "weakening" arises solely from testing (projecting) the strong residual against basis functions and penalizing the resulting moments. This choice is deliberate: it yields a general, operator-agnostic construction applicable to linear or nonlinear PDEs without requiring case-by-case derivations of integrated-by-parts bilinear forms and associated boundary flux terms. Since $u_N$ is represented in a smooth spectral basis and boundary conditions are enforced by construction (e.g., via boundary-satisfying bases or lifting), the strong residual and its moments are well-defined and can be differentiated efficiently via automatic differentiation.

To incorporate the initial condition $u(x, y, 0) = u_0(x, y)$, VSL adds a separate initial-condition term

$$\mathcal{L}_{\text{IC}}(\mathbf{c}) = \frac{1}{|\mathcal{S}_{\text{IC}}|} \sum_{(x_i, y_i) \in \mathcal{S}_{\text{IC}}} \left( u_N(x_i, y_i, 0; \mathbf{c}) - u_0(x_i, y_i) \right)^2, \tag{85}$$

where $\mathcal{S}_{\text{IC}}$ is a set of spatial collocation points at $t = 0$, typically chosen as Chebyshev–Gauss interior nodes in each coordinate. The total VSL loss has the form

$$\mathcal{J}_{\text{VSL}}(\mathbf{c}) = E_{\text{strong/weak(GLS)}}(\mathbf{c}) + \lambda_{\text{IC}} \mathcal{L}_{\text{IC}}(\mathbf{c}) + \lambda_{\text{reg}} \|\mathbf{c}\|_2^2, \tag{86}$$

where $\lambda_{\text{IC}}$ and $\lambda_{\text{reg}}$ are weights on the initial-condition and regularization terms, respectively. Coefficients are initialized, the loss is minimized using an adaptive optimizer such as Adam, and training dynamics are shaped by a cosine decay restarts learning-rate schedule.

To provide a classical high-order reference, the benchmark includes a spectral collocation solver that discretizes space with Chebyshev–Lobatto points and advances in time using the Crank–Nicolson scheme. In each spatial direction, say $x$, one constructs Lobatto nodes $z_j = \cos(\pi j / N_x^{\text{spec}})$ on $[-1, 1]$ and maps them to $[0, 1]$ via $x_j = (z_j + 1)/2$. A Chebyshev differentiation matrix $D_x$ is built using standard formulas involving the scaling vector and pairwise differences of the nodes; the second-derivative matrix $D_{xx}$ in the physical coordinate is obtained as $D_{xx} = 4 D_x^2$ to account for the affine map between $[-1, 1]$ and $[0, 1]$. The same construction is applied in $y$, yielding $D_{yy}$. Restricting to interior indices and forming Kronecker products with identity matrices produces a discrete Laplacian $L$ acting on interior degrees of freedom. The semi-discrete system for interior values $u_{\text{int}}(t)$ can be written as

$$\frac{d}{dt} u_{\text{int}}(t) = \nu L u_{\text{int}}(t) + f_{\text{int}}(t), \tag{87}$$

where $f_{\text{int}}(t)$ are the interior values of the forcing term. On a uniform time grid $\{t_n\}_{n=0}^{N_t^{\text{steps}}}$ with step size $\Delta t$, $I$ is the interior identity matrix, Crank–Nicolson updates are given by

$$\left(I - \frac{\Delta t}{2}\nu L\right) u_{\text{int}}^{n+1} = \left(I + \frac{\Delta t}{2}\nu L\right) u_{\text{int}}^n + \frac{\Delta t}{2}(f_{\text{int}}^{n+1} + f_{\text{int}}^n), \tag{88}$$

with initial data $u_{\text{int}}^0$ obtained by evaluating $u_0(x, y)$ at interior collocation nodes. The boundary values are set to zero at each time step to enforce the homogeneous Dirichlet conditions. This scheme is second-order accurate in time, high-order accurate in space, and unconditionally stable for the diffusion equation, thus providing a stringent baseline against which VSL and PINN approaches can be compared.

The PINN baseline for this benchmark adopts a field-space neural parameterization with analytic lifting to enforce both boundary and initial conditions. A neural network $v_\theta(x, y, t)$ with several hidden layers and smooth activation functions (e.g., tanh) is defined on $\mathbb{R}^3$. The physical solution is represented as

$$u_\theta(x, y, t) = (1 - t)u_0(x, y) + t\, x(1 - x)y(1 - y)\, v_\theta(x, y, t), \tag{89}$$

where $u_0(x, y) = \sin(\pi x)\sin(\pi y)$ is the prescribed initial profile. This lifting ensures that at $t = 0$,

$$u_\theta(x, y, 0) = u_0(x, y), \tag{90}$$

and along the spatial boundary $x \in \{0,1\}$ or $y \in \{0,1\}$, the factor $x(1-x)y(1-y)$ vanishes, giving $u_\theta = 0$ for all $t$. Therefore, the PDE constraints are the only ones that require explicit enforcement in the loss. The PINN loss is defined as a mean-squared PDE residual evaluated at collocation points $\{(x_i, y_i, t_i)\}$ in space–time,

$$\mathcal{L}_{\text{PINN}}(\theta) = \frac{1}{|\mathcal{S}_{\text{pde}}|} \sum_{(x_i, y_i, t_i) \in \mathcal{S}_{\text{pde}}} \left( u_{\theta,t}(x_i, y_i, t_i) - \nu\left(u_{\theta,xx}(x_i, y_i, t_i) + u_{\theta,yy}(x_i, y_i, t_i)\right) - f(x_i, y_i, t_i) \right)^2, \tag{91}$$

with derivatives again obtained via automatic differentiation. Although the initial condition is enforced analytically, a separate initial-condition diagnostic loss may be monitored to quantify numerical deviations due to autodifferentiation and implementation details. The PINN is trained with the same optimizer and cosine decay restarts learning-rate schedule used in VSL, to ensure that differences in performance are attributable to the parameterization and loss structure rather than to optimization choices.

The discussion of the 2D diffusion benchmark emphasizes the conceptual contrasts between the three approaches. VSL treats the diffusion problem as a global space–time variational problem in a carefully constructed spectral subspace, with boundary conditions satisfied by construction and initial conditions enforced by a dedicated loss term. The space–time perspective has the advantage that no separate time-stepping procedure is required: the solution is approximated as a single function over $Q$, and all temporal dynamics are captured through the spectral time basis and the space–time energy. This contrasts with the classical collocation baseline, which employs a semi-discrete strategy: it discretizes space spectrally and uses Crank–Nicolson to advance in time, effectively splitting the problem into a sequence of elliptic solves. While the classical approach is highly mature and efficient, it is less naturally suited to scenarios in which time must be treated on an equal footing with other coordinates, such as parametric or stochastic dimensions. The PINN baseline, in turn, offers great flexibility in principle, as the neural architecture can, in theory, adapt to complex solution structures and be extended to more general geometries. However, its optimization landscape is more complex because the parameterization is high-dimensional and nonlinear, and the PDE residual loss generally needs to balance multiple contributions and scales. In contrast, VSL's linear coefficient-space parameterization and physics-derived energies lead to an optimization problem that is closer to a smooth nonlinear least-squares system with a fixed feature map, aligning more closely with classical numerical analysis.

We start our numerical investigation with the 1D heat equation with $\nu = 1$. We consider the manufactured 1D heat equation on $(0,1) \times (0,1)$ with exact solution $u^\star(x, t) = e^{-t}\sin(\pi x)$ and matching forcing. Three solvers are compared: (i) VSL in a tensor–product Chebyshev space with $(N_x = N_t = 8)$, trained by Adam with cosine–decay restarts; (ii) a classical reference method based on Chebyshev–Lobatto collocation in space and Crank–Nicolson in time with $(N_x^{\text{spec}} = 32)$ and $(N_t = 64)$; and (iii) a PINN with lifting that enforces the Dirichlet boundary conditions and the initial condition analytically.

The classical collocation method (see figures 9-12) achieves the expected spectral-in-space/high-order accuracy. At $t = 1$, it yields relative errors $L_{\text{rel}}^2 = 2.284 \times 10^{-6}$ and $L_{\text{rel}}^\infty = 2.316 \times 10^{-6}$. This is consistent with a well-resolved Chebyshev–Crank–Nicolson reference on this smooth problem and confirms that the baseline is numerically reliable.

In the GLS weak setting (see figures 9 and 10), the optimization minimizes the squared moments of the strong residual against the tensor–product trial/test basis (Galerkin projection energy). Training reduces the GLS weak energy to very small values (e.g., $E_{\text{weak(GLS)}} \sim 10^{-5} - 10^{-8}$ in later epochs) while the diagnostic strong residual oscillates due to the restart schedule. At final time, VSL achieves $L_{\text{rel}}^2 = 7.929 \times 10^{-5}$ and $L_{\text{rel}}^\infty = 1.116 \times 10^{-4}$. Switching to the strong-form VSL objective (see figures 11 and 12) leads to improved solution quality with the same spectral rank and optimizer settings. At $t = 1$, we obtain $L_{\text{rel}}^2 = 5.646 \times 10^{-5}$ and $L_{\text{rel}}^\infty = 8.299 \times 10^{-5}$, and the space–time max error drops sharply to $\max(|u_{\text{VSL}} - u^\star|) \approx 5.86 \times 10^{-5}$.

The PINN training (see figures 9-12) uses a lifted ansatz so that the initial condition and boundary conditions are satisfied by construction, which is reflected in the reported IC diagnostic being identically zero. Nevertheless, the PDE-residual loss decreases more slowly than for VSL, and the final-time errors remain larger than both VSL variants: in the GLS weak experiment $L^2_{\text{rel}} = 4.860 \times 10^{-4}$, $L^\infty_{\text{rel}} = 5.043 \times 10^{-4}$, and $\max(|u_{\text{PINN}} - u^\star|) \approx 7.96 \times 10^{-4}$; in the strong experiment $L^2_{\text{rel}} = 4.133 \times 10^{-4}$, $L^\infty_{\text{rel}} = 3.883 \times 10^{-4}$, and $\max(|u_{\text{PINN}} - u^\star|) \approx 7.55 \times 10^{-4}$. This is consistent with the well-known stiffness of PINN optimization on even simple parabolic PDEs and highlights the advantage of solving directly in a structured spectral coefficient space.

For 2D heat equation with $\nu = 0.1$, figures 13 and 14 summarize the performance of the coefficient-space VSL model trained with the space–time GLS weak-form energy using a compact tensor-product Chebyshev basis $(N_x, N_y, N_t) = (6,6,6)$. At $t = 1$, the GLS weak-form solution reproduces the correct spatial structure of the manufactured solution $u^\star(x,y,t) = e^{-t}\sin(\pi x)\sin(\pi y)$ and yields relative errors $L^2_{\text{rel}} = 1.028 \times 10^{-3}$ and $L^\infty_{\text{rel}} = 1.807 \times 10^{-3}$. The classical Chebyshev–Lobatto collocation + Crank–Nicolson solver provides a much tighter reference with $1.30 \times 10^{-5}$ in both relative $L^2_{\text{rel}}$ and $L^\infty_{\text{rel}}$, and a very small maximum absolute error $\approx 4.78 \times 10^{-6}$, while the PINN (width $(= 64)$, depth $(= 4)$) attains $2.56 \times 10^{-3}$ ($L^2_{\text{rel}}$) and $5.041 \times 10^{-3}$ ($L^\infty_{\text{rel}}$). Thus, in this configuration, GLS weak-form VSL is clearly more accurate than the PINN at the final time, though it remains about two orders of magnitude less accurate than the classical spectral reference.

The diagnostics in figure 14 clarify the GLS weak-form training behavior. The weak-form energy $E_{\text{weak(GLS)}}$ decreases by multiple orders of magnitude and reaches extremely small values (down to $\sim 10^{-9}$ around epoch $\approx 1400$), indicating that the Galerkin moments of the residual have been strongly suppressed. At the same time, the pointwise strong residual diagnostic (mean-square PDE residual evaluated at Chebyshev interior collocation points) is not monotonically decreasing early in training and can even be large while the GLS weak energy is already small (e.g., the diagnostic PDE residual peaks near epoch 100). This decoupling is expected: the GLS weak objective drives residual orthogonality to the chosen test/basis functions, whereas the diagnostic monitors pointwise residual magnitude on a separate set. Overall, the GLS weak-form VSL objective is effective at enforcing the PDE in an integrated (Galerkin) sense and produces a high-quality solution with modest basis size, but the final accuracy is ultimately limited by the combination of the GLS weak objective and the low-dimensional coefficient representation used here.

Switching VSL to the strong-form space–time residual least-squares objective produces a marked improvement in accuracy for the same basis resolution (6,6,6). At $t = 1$, strong-form VSL attains $L^2_{\text{rel}} = 4.806 \times 10^{-5}$ and $L^\infty_{\text{rel}} = 7.539 \times 10^{-5}$, substantially better than the GLS weak-form run and much closer to the classical Chebyshev–Lobatto collocation + Crank–Nicolson reference ($1.30 \times 10^{-5}$ in both norms). In contrast, the PINN performance remains essentially unchanged between the two campaigns ($L^2_{\text{rel}} \approx 2.52 \times 10^{-3}$, $L^\infty_{\text{rel}} \approx 4.81 \times 10^{-3}$), consistent with the fact that only the VSL objective was altered. The solution and error fields in figure 15 reflect this: strong-form VSL significantly reduces the final-time error magnitude compared with GLS weak-form, while preserving the correct qualitative structure.

The training traces in figure 16 show that strong-form VSL aligns the optimized objective more directly with the diagnostic being monitored. The strong energy $E_{\text{strong}}$ decreases rapidly (from $1.025 \times 10^{-1}$ at epoch 0 to $\sim 10^{-8}$ around epoch $\approx 1400$), and the collocation-based PDE residual diagnostic falls in tandem to $\sim 10^{-7}$ at its best before the learning-rate restart. Overall, the strong-form objective provides a more direct route to minimizing the pointwise residual and delivers near-classical accuracy at $t = 1$ with a compact coefficient model.

Finally, Benchmark II demonstrates how VSL can act as a bridge between fully variational, space–time formulations and modern learning infrastructure. It shows that a global spectral representation in space–time, combined with strong or GLS weak variational energies and gradient-based optimization, can reproduce the behavior of classical solvers while retaining compatibility with the automatic differentiation and scheduling tools of deep learning. At the same time, comparisons with PINNs highlight the trade-offs between structured, interpretable coefficient-space learning and more flexible, but less analyzable, field-space neural solvers. This benchmark therefore plays a central role in elucidating the strengths and limitations of VSL for time-dependent PDEs and sets the stage for more challenging scenarios such as nonlinear diffusion–advection systems considered in subsequent benchmarks.

## 8. Benchmark III – Burgers-Type Diffusion/Advection Problems

Benchmark III addresses Burgers-type diffusion–advection problems and is designed to probe the behavior of VSL in the presence of nonlinear transport terms superimposed on diffusive dynamics [58]. In contrast to the purely elliptic and linear parabolic model problems in the first two benchmarks, Burgers-type equations introduce a convective nonlinearity that couples the solution to its own gradient. Even for smooth manufactured solutions, this nonlinearity leads to steeper gradients, potentially more intricate residual landscapes, and nonlinear algebraic systems whose Jacobians depend on both the solution and its derivatives. The benchmark therefore examines whether the coefficient-space VSL formulation, which was previously tested on linear operators, remains stable and accurate when the PDE operator becomes nonlinear, and how its

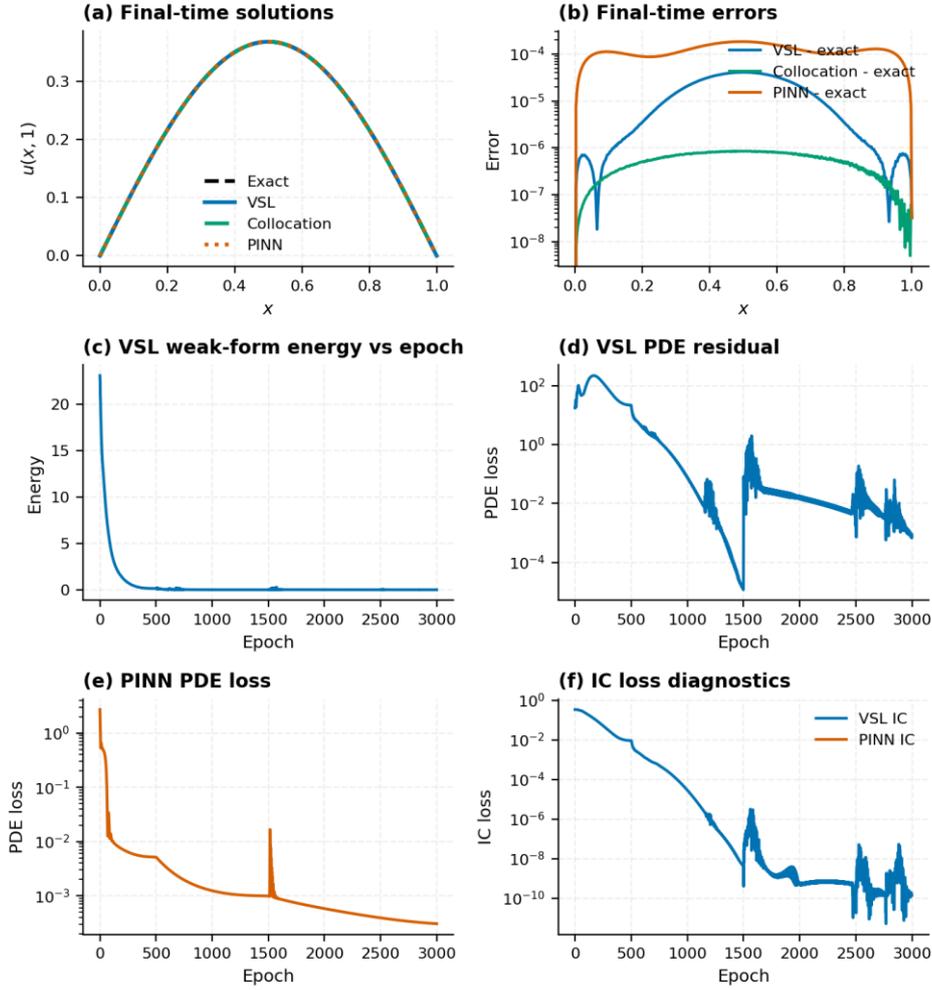

**Figure 9.** Overview of final-time accuracy and training diagnostics for the 1D heat-equation benchmark. (a) Final-time profiles $u(x,1)$: exact solution, VSL, classical collocation, and PINN. (b) Final-time absolute errors $|u(x,1) - u^\star(x,1)|$ shown on a logarithmic scale. (c) VSL objective history (weak moment-residual energy $E_{\text{weak}}$). (d) VSL diagnostic strong-form PDE residual loss evaluated on a fixed collocation set. (e) PINN PDE residual loss during training. (f) Initial-condition diagnostic losses for VSL and PINN.

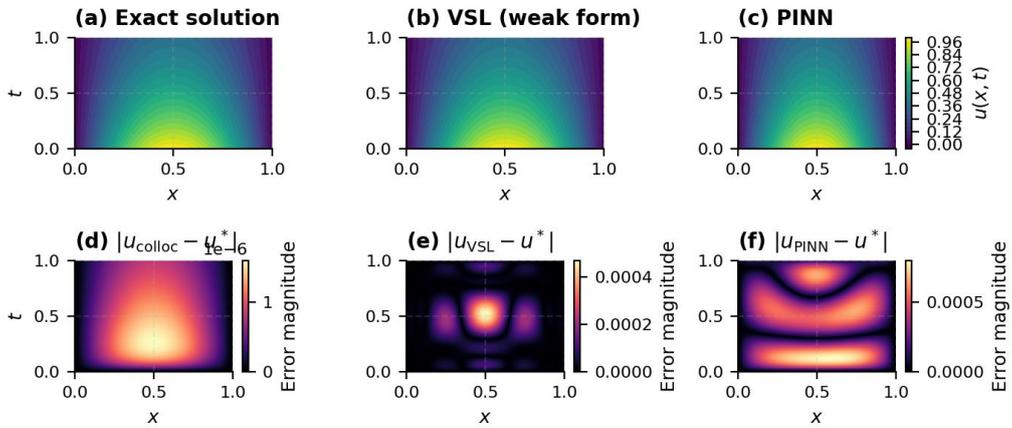

**Figure 10.** Space–time comparison for the 1D heat equation on $(x,t) \in [0,1] \times [0,1]$. Top row: (a) exact solution $u^\star(x,t) = e^{-t}\sin(\pi x)$; (b) VSL prediction $u_{\text{VSL}}(x,t)$ (moment-residual/Galerkin "weak" objective); (c) PINN prediction $u_{\text{PINN}}(x,t)$ with analytic IC/BC lifting. Bottom row: absolute error fields (d) $|u_{\text{colloc}} - u^\star|$ for Chebyshev–Lobatto collocation in $x$ with Crank–Nicolson time stepping; (e) $|u_{\text{VSL}} - u^\star|$; (f) $|u_{\text{PINN}} - u^\star|$.

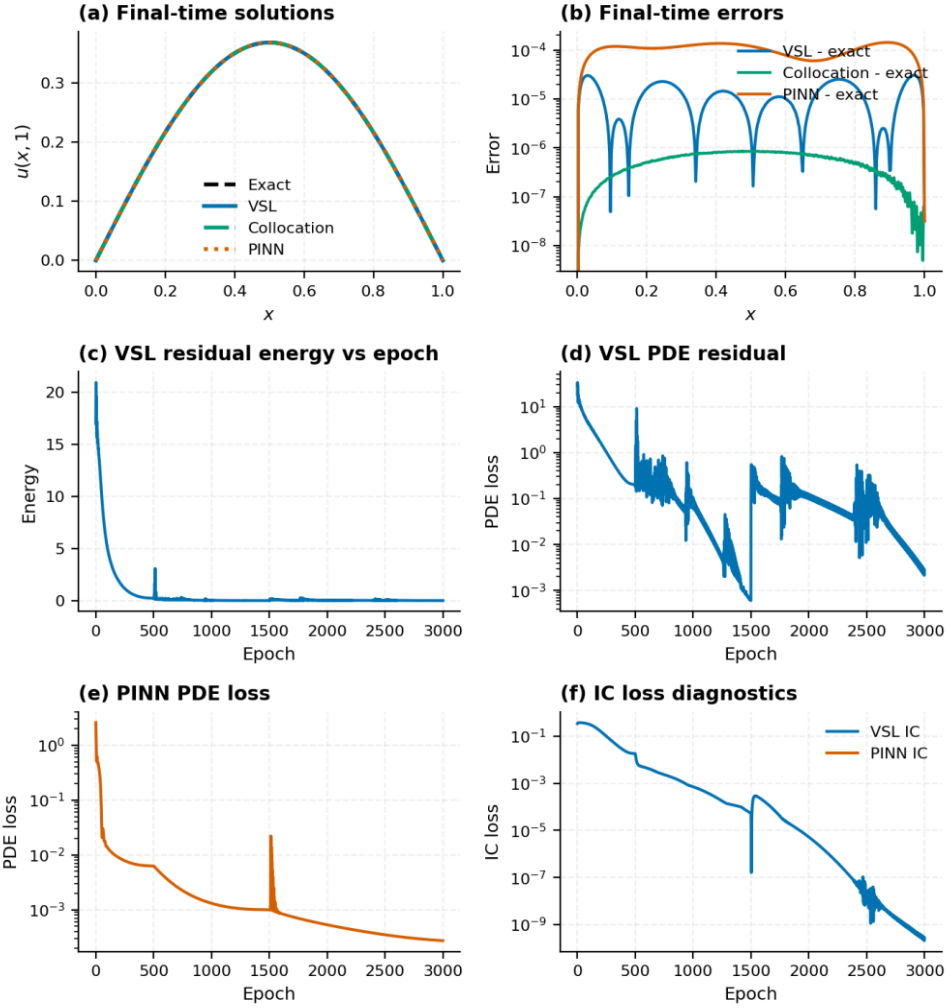

**Figure 11.** Space–time comparison for the 1D heat equation on $(x,t) \in [0,1] \times [0,1]$. The same as in figure 9, but for VSL strong-form objective.

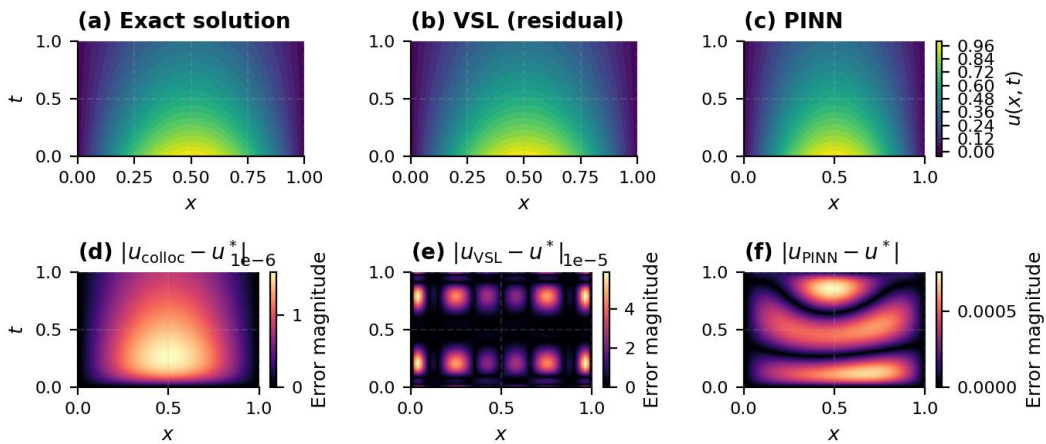

**Figure 12.** Space–time comparison for the 1D heat equation on $(x,t) \in [0,1] \times [0,1]$. The same as in figure 10, but for VSL strong-form objective.

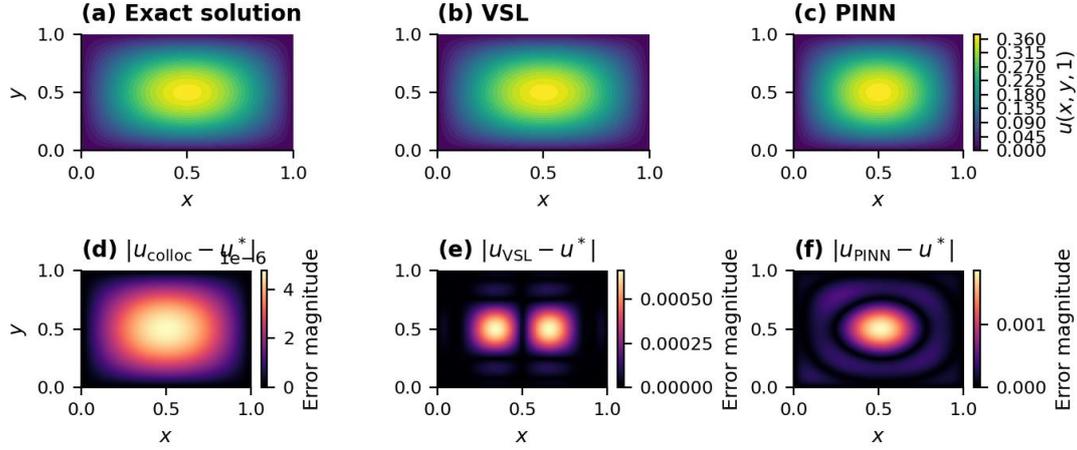

**Figure 13.** (Weak-form VSL) Final-time solutions and error fields for the 2D time-dependent diffusion benchmark at $t = 1$. (a–c) Solution fields on a common $N_x \times N_y$ test grid: (a) manufactured exact solution $u^\star(x, y, 1) = e^{-1} \sin(\pi x) \sin(\pi y)$, (b) VSL prediction $u_{\text{VSL}}(x, y, 1)$ obtained by minimizing a space–time Galerkin (weak-form) least-squares energy, and (c) PINN prediction $u_{\text{PINN}}(x, y, 1)$ with analytic lifting for BC/IC. (d–f) Absolute error heatmaps at $t = 1$: (d) classical Chebyshev–Lobatto collocation + Crank–Nicolson, (e) VSL, and (f) PINN, each compared against $u^\star$.

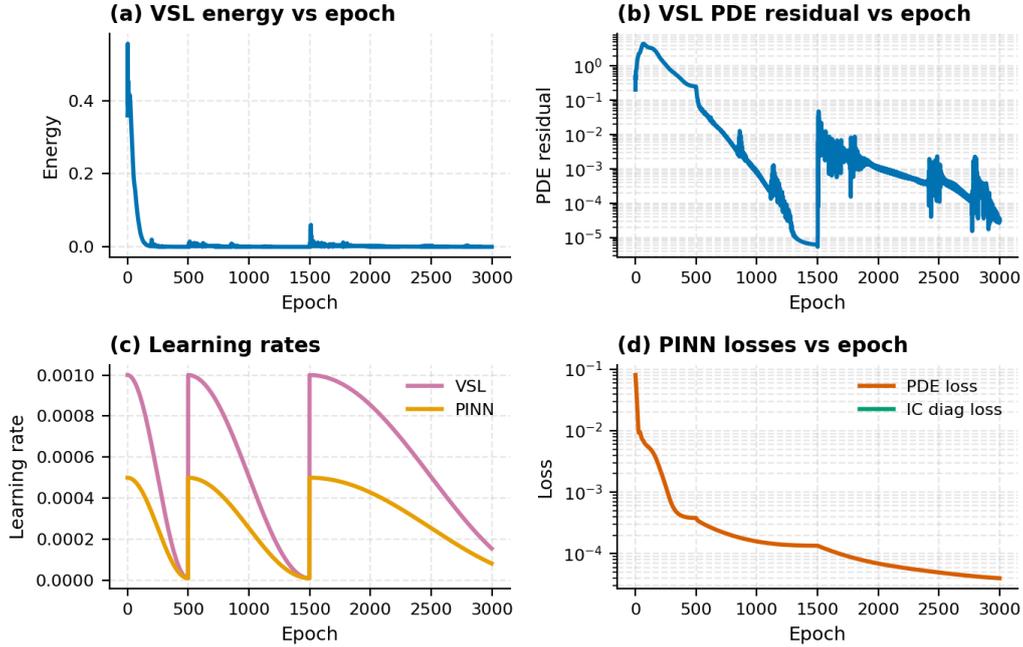

**Figure 14.** (Weak-form VSL) Training diagnostics for the 2D time-dependent diffusion campaign (VSL vs PINN). (a) VSL training energy versus epoch, where the energy corresponds to the weak-form space–time Galerkin least-squares objective used in coefficient-space training. (b) VSL diagnostic PDE residual loss at Chebyshev interior collocation points (mean-square of $u_t - \nu \Delta u - f$), tracked during training. (c) CosineDecayRestarts learning-rate schedules for VSL and PINN optimizers. (d) PINN training losses: PDE residual loss and initial-condition diagnostic loss (both in log scale), where BC and IC are enforced analytically by lifting.

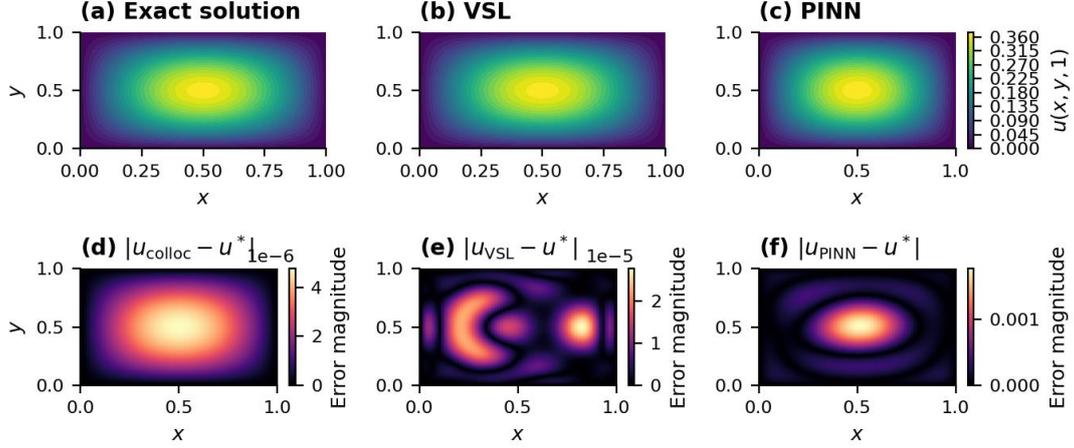

**Figure 15.** The same as in figure 13, but for VSL strong-form objective.

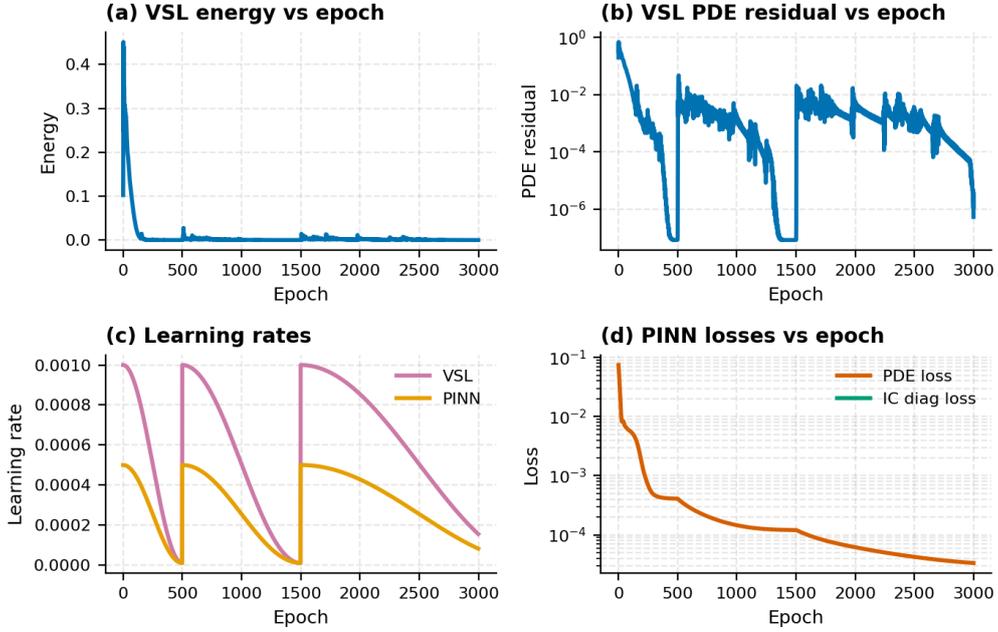

**Figure 16.** The same as in figure 14, but for VSL strong-form objective.

behavior compares to a classical Chebyshev collocation solver and PINN trained directly in physical space.

The one-dimensional (1D) component of Benchmark III considers the stationary Burgers-type boundary-value problem on the unit interval [58],

$$\nu u_{xx}(x) - u(x)u_x(x) = f(x), \qquad x \in (0,1), u(0) = u(1) = 0, \tag{92}$$

where $\nu > 0$ is a viscosity parameter and $f$ is chosen by the method of manufactured solutions. Specifically, the exact solution is prescribed as

$$u^\star(x) = \sin(\pi x), \tag{93}$$

and the forcing is defined so that $u^\star$ satisfies the strong PDE and the homogeneous Dirichlet boundary conditions exactly. With

$$u^\star(x) = \sin(\pi x), \qquad u^\star_x(x) = \pi \cos(\pi x), \qquad u^\star_{xx}(x) = -\pi^2 \sin(\pi x), \tag{94}$$

the forcing term becomes

$$f(x) = \nu u^\star_{xx}(x) - u^\star(x)u^\star_x(x) = -\nu\pi^2 \sin(\pi x) - \pi \sin(\pi x)\cos(\pi x). \tag{95}$$

This manufactured solution construction guarantees that the exact solution is known in closed form, enabling direct computation of both $L^2_{\text{rel}}$ and $L^\infty_{\text{rel}}$ errors for all numerical methods on a dense test grid [59].

The two-dimensional (2D) part of the benchmark extends the same strategy to the steady Burgers-type diffusion–advection equation on the unit square $\Omega = (0,1)^2$,

$$\nu(u_{xx} + u_{yy})(x,y) - u(x,y)u_x(x,y) - u(x,y)u_y(x,y) = f(x,y), \qquad (x,y) \in \Omega, \tag{96}$$

supplemented with homogeneous Dirichlet boundary conditions

$$u(x,y) = 0, \qquad (x,y) \in \partial\Omega. \tag{97}$$

Here the reference solution is chosen as

$$u^\star(x,y) = \sin(\pi x)\sin(\pi y), \tag{98}$$

from which the forcing is derived by direct substitution into the PDE. one obtains

$$\begin{aligned}f(x,y) &= \nu(u^\star_{xx} + u^\star_{yy}) - u^\star u^\star_x - u^\star u^\star_y \\ &= -2\nu\pi^2 \sin(\pi x)\sin(\pi y) - \pi\sin(\pi x)\sin(\pi y)(\cos(\pi x)\sin(\pi y) + \sin(\pi x)\cos(\pi y)).\end{aligned} \tag{99}$$

As in the 1D case, this setup yields an analytically known solution and provides a controlled environment for comparing strong- and GLS weak-form VSL, classical tensor-product Chebyshev collocation, and a 2D PINN.

In the present implementation, the VSL energies are built from the strong residual evaluated at quadrature nodes, and the connection to the GLS weak form is realized through weighted residuals and Galerkin projections rather than by assembling a traditional stiffness matrix. Within the VSL framework, both the 1D and 2D problems are parameterized in coefficient space using Dirichlet-satisfying Chebyshev bases. In 1D, the approximate solution is expanded as

$$u_{N_x}(x) = \sum_{k=0}^{N_x-1} c_k \phi_k(x), \tag{100}$$

where the basis functions $\phi_k(x) = T_{k+2}(2x-1) - T_k(2x-1)$, $k = 0, \ldots, N_x - 1$, are built from Chebyshev polynomials of the first kind $T_n$. By construction, $\phi_k(0) = \phi_k(1) = 0$ for all $k$, so any coefficient vector $\mathbf{c}$ produces a field $u_{N_x}$ that exactly satisfies the homogeneous Dirichlet boundary conditions. In 2D, a tensor-product basis is formed from identical 1D Dirichlet Chebyshev bases in $x$ and $y$,

$$\Phi_{k,\ell}(x,y) = \phi_k^{(x)}(x)\phi_\ell^{(y)}(y), \qquad k = 0, \ldots, N_x - 1;\ \ell = 0, \ldots, N_y - 1, \tag{101}$$

yielding a spectral approximation

$$u_{N_x,N_y}(x,y) = \sum_{k=0}^{N_x-1} \sum_{\ell=0}^{N_y-1} c_{k\ell} \Phi_{k,\ell}(x,y). \tag{102}$$

The boundary conditions $u = 0$ on $\partial\Omega$ are therefore enforced analytically by the choice of basis, without the need for explicit penalty terms or boundary losses.

Given this spectral representation, VSL defines a strong residual of the form

$$r(x;\mathbf{c}) = \nu u_{N_x,xx}(x;\mathbf{c}) - u_{N_x}(x;\mathbf{c})u_{N_x,x}(x;\mathbf{c}) - f(x), \tag{103}$$

in 1D and

$$r(x,y;\mathbf{c}) = \nu\left(u_{N_x,N_y,xx} + u_{N_x,N_y,yy}\right)(x,y;\mathbf{c}) - u_{N_x,N_y}(x,y;\mathbf{c})\left(u_{N_x,N_y,x} + u_{N_x,N_y,y}\right)(x,y;\mathbf{c}) - f(x,y), \tag{104}$$

in 2D. These residuals are evaluated at Gauss–Legendre quadrature nodes, and the required derivatives are obtained by automatic differentiation using nested tape constructs. Two variational energies are then considered. The strong-form VSL energy minimizes the squared $L^2$ norm of the residual,

$$E_{\text{strong}}(\mathbf{c}) = \int_0^1 (r(x;\mathbf{c}))^2 dx, \qquad E_{\text{strong}}(\mathbf{c}) = \iint_\Omega (r(x,y;\mathbf{c}))^2 dxdy, \tag{105}$$

approximated by Gauss–Legendre quadrature as

$$E_{\text{strong}}(\mathbf{c}) \approx \sum_{q_x=1}^{Q_x} w_{q_x}\left(r(x_{q_x};\mathbf{c})\right)^2, \qquad E_{\text{strong}}(\mathbf{c}) \approx \sum_{q_x=1}^{Q_x} \sum_{q_y=1}^{Q_y} w_{q_x} w_{q_y}\left(r\left(x_{q_x}, y_{q_y};\mathbf{c}\right)\right)^2. \tag{106}$$

In parallel, a GLS weak VSL energy is defined by projecting the strong residual onto the spectral basis and minimizing the squared norm of the resulting coefficient vector. In 1D, defining

$$R_k(\mathbf{c}) = \int_0^1 r(x;\mathbf{c})\phi_k(x)dx, \tag{107}$$

the weak energy reads

$$E_{\text{weak(GLS)}}(\mathbf{c}) = \frac{1}{2} \sum_{k=0}^{N_x-1} (R_k(\mathbf{c}))^2, \tag{108}$$

with integrals again approximated by quadrature. In 2D, the Galerkin residuals

$$R_{k,\ell}(\mathbf{c}) = \iint_\Omega r(x,y;\mathbf{c})\Phi_{k,\ell}(x,y)dxdy, \tag{109}$$

are assembled at tensor-product Gauss–Legendre nodes, and the weak energy is

$$E_{\text{weak(GLS)}}(\mathbf{c}) = \frac{1}{2}\sum_{k=0}^{N_x-1}\sum_{\ell=0}^{N_y-1}(R_{k,\ell}(\mathbf{c}))^2. \tag{110}$$

In both strong and GLS weak modes, the coefficient vector $\mathbf{c}$ is treated as a trainable parameter vector and is optimized using Adam with a cosine-decay-restarts learning-rate schedule. A separate diagnostic loss monitors the mean-squared strong residual evaluated at Chebyshev interior collocation nodes, providing an independent measure of how well the learned coefficients satisfy the PDE away from the quadrature points used in the energy.

The classical reference solver used for comparison in 1D is a Chebyshev–Lobatto collocation method with a damped Newton iteration. Chebyshev–Lobatto nodes $z_k = \cos(\pi k/N)$ on $[-1,1]$ are mapped to physical coordinates $x = (z+1)/2$, and standard Trefethen-style differentiation matrices are assembled to approximate $d/dx$ and $d^2/dx^2$. The unknowns are the interior nodal values $u(x_j)$ at $j = 1, \ldots, N-1$, while the boundary values $u(0) = u(1) = 0$ are imposed exactly. The discrete residual at each node enforces the strong PDE,

$$r_j = \nu(u_{xx})_j - u_j(u_x)_j - f_j, \tag{111}$$

and the Jacobian of this nonlinear system is constructed analytically using the product rule and the structure of the differentiation matrices. A damped Newton method with a simple backtracking line search is used to solve for the interior degrees of freedom. Once converged, the collocation solution is interpolated to a dense test grid using barycentric Chebyshev interpolation for error evaluation and plotting. In 2D, the same approach is extended to a tensor-product Chebyshev–Lobatto grid in $x$ and $y$. Differentiation matrices in each direction are combined to approximate $u_x, u_y, u_{xx}, u_{yy}$, and the nonlinear PDE is enforced at interior grid points with $u = 0$ imposed on all boundaries. The resulting large but structured nonlinear system is again solved with damped Newton, and the converged grid solution is interpolated to a uniform test grid for comparison with VSL and PINN.

The PINN baselines in Benchmark III employ multilayer perceptron trained to minimize the mean-squared strong-form residual in physical space, with homogeneous Dirichlet boundary conditions enforced analytically by boundary lifting. In 1D, a scalar network $v_\theta(x)$ is constructed, and the physical solution is represented as

$$u(x) = x(1-x)v_\theta(x), \tag{112}$$

which guarantees $u(0) = u(1) = 0$ for all network parameters. The PDE loss is defined as

$$\mathcal{L}_{\text{PINN}}(\theta) = \mathbb{E}_x\left[\left(\nu u_{xx}(x) - u(x)u_x(x) - f(x)\right)^2\right], \tag{113}$$

with derivatives computed via automatic differentiation and the expectation approximated by Chebyshev interior collocation nodes. An additional diagnostic boundary loss can be monitored, but for the lifted representation it remains negligible. In 2D, a network $v_\theta(x,y)$ with two inputs is used, and the lifted representation

$$u(x,y) = x(1-x)y(1-y)v_\theta(x,y), \tag{114}$$

enforces the homogeneous Dirichlet condition on all four edges. The PINN loss in 2D is

$$\mathcal{L}_{\text{PINN}}(\theta) = \mathbb{E}_{(x,y)}\left[\left(\nu(u_{xx} + u_{yy}) - uu_x - uu_y - f(x,y)\right)^2\right], \tag{115}$$

again approximated at Chebyshev interior collocation points. In both dimensions, the networks are trained with Adam and cosine-decay restarts, and loss histories are recorded for analysis.

Comparisons across the three families of solvers are carried out using both quantitative and qualitative diagnostics. Quantitatively, all methods are evaluated against the manufactured exact solution on dense uniform grids (400 points in 1D and ($64 \times 64$) in 2D). $L^2_{\text{rel}}$ and $L^\infty_{\text{rel}}$ norms are reported for VSL, classical collocation, and PINN. Qualitatively, the benchmark uses a suite of diagnostic plots: overlayed solution profiles in 1D, log-scale pointwise absolute error curves, VSL energy and strong residual histories, PINN loss histories, and spectra of the learned Chebyshev coefficients $c_k$. In 2D, contour and heatmap visualizations of the solutions and error fields, as well as midline slices at fixed $y$, allow a direct visual assessment of how accurately each method reproduces the structure of the exact solution across the domain.

We start our numerical investigation with the 1D Burgers equation (see figures 17 and 18). We benchmarked three solvers for the stationary 1D Burgers-type boundary-value problem with manufactured solution $u_{\text{exact}}(x) = \sin(\pi x)$: (i) VSL in coefficient space using a Dirichlet-satisfying Chebyshev basis with $N_x = 16$ modes, (ii) classical Chebyshev–Lobatto collocation with a damped Newton solver on $N_{\text{spec}} = 32$ nodes, and (iii) a boundary-lifted PINN (Keras MLP with width 64 and depth 4). In the GLS weak-form VSL variant, the optimization minimizes the least-squares Galerkin residuals. Training logs indicate that the GLS weak-form energy decreased substantially from $E \approx 1.836$ at initialization to $E \approx 3.05 \times 10^{-4}$ by epoch 2900, with the diagnostic strong residual loss dropping from 5.859 initially to $2.253 \times 10^{-4}$ at epoch 2900. Quantitatively, GLS weak form achieved relative errors $L^2_{\text{rel}} = 7.605 \times 10^{-5}$ and $L^\infty_{\text{rel}} = 1.513 \times 10^{-4}$. This accuracy is strong for a compact spectral representation $N_x = 16$ and confirms that minimizing projected (Galerkin) residuals is sufficient to recover an accurate solution even when the diagnostic strong residual is not enforced directly. The classical collocation method provided the most accurate and fastest baseline: $L^2_{\text{rel}} = 7.653 \times 10^{-8}$, and $L^\infty_{\text{rel}} = 2.038 \times 10^{-7}$. This is consistent with the spectral accuracy expected for a smooth manufactured solution when a Newton solver converges. The PINN baseline converged smoothly to a PDE loss $\approx 6.66 \times 10^{-6}$ by epoch 2900 and achieved $L^2_{\text{rel}} = 4.186 \times 10^{-4}$, $L^\infty_{\text{rel}} = 5.296 \times 10^{-4}$. Boundary lifting ensured the BC diagnostic loss remained essentially zero throughout, isolating the learning burden to the interior PDE residual. Overall, the GLS weak-form VSL results demonstrate a clear advantage over the PINN in accuracy for this configuration while remaining far behind classical collocation in both accuracy and runtime.

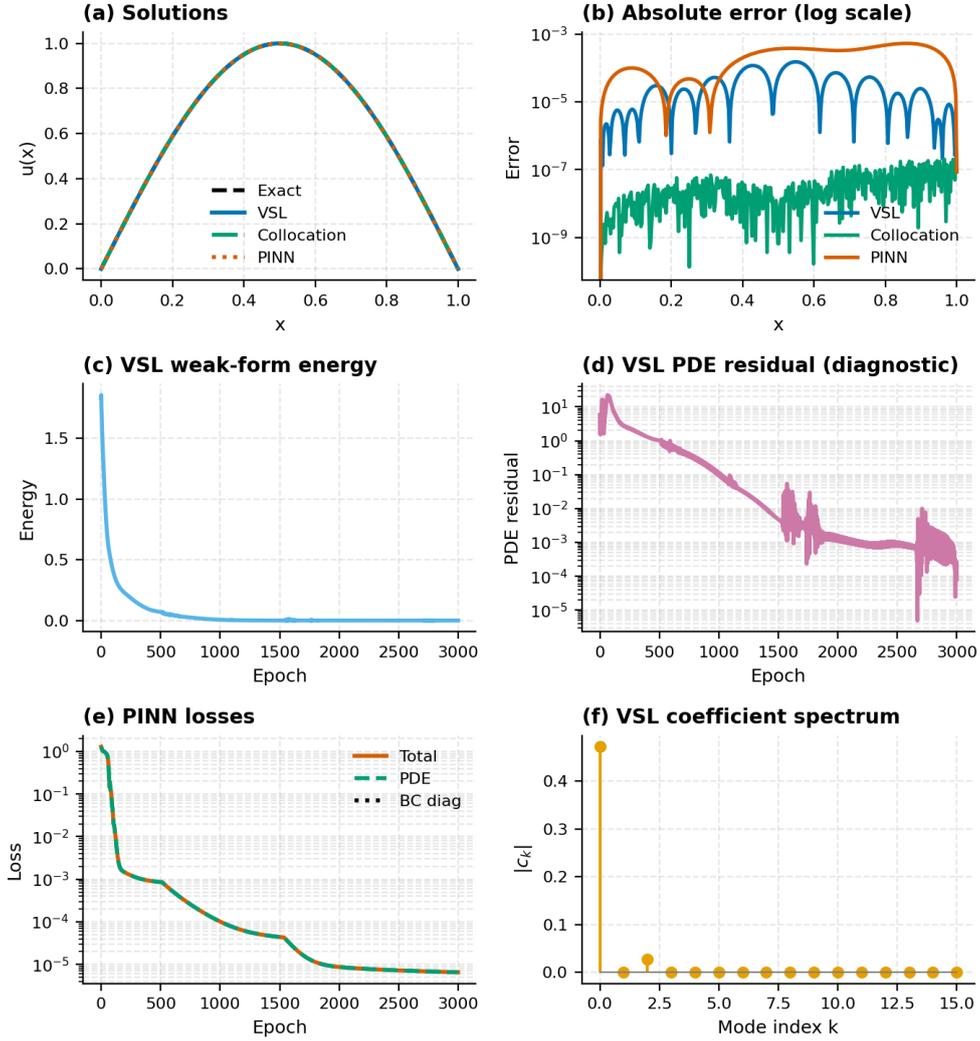

**Figure 17.** Overview of the 1D Burgers-type boundary-value benchmark comparing VSL, classical Chebyshev collocation, and PINN. (a) Numerical solutions on a dense grid (400 points): exact manufactured solution $u(x) = \sin(\pi x)$ versus VSL (Dirichlet-satisfying Chebyshev basis, $N_x = 16$, classical Chebyshev–Lobatto collocation $N_{\text{spec}} = 32$, and a boundary-lifted PINN (MLP with width 64 and depth 4). (b) Pointwise absolute errors $|u_{\text{num}} - u_{\text{exact}}|$ on a logarithmic scale for the three methods. (c) VSL training curve for the GLS weak-form objective $E_{\text{weak(GLS)}}$ (**c**). (d) Diagnostic VSL strong-form residual loss evaluated at interior Chebyshev points during training (log scale). (e) PINN training losses (total loss, PDE residual loss, and a near-zero boundary-condition diagnostic due to analytical boundary lifting $u(x) = x(1-x)v(x)$. (f) Magnitude spectrum of the learned VSL Chebyshev coefficients $|c_k|$, illustrating spectral decay.

We repeated the same Burgers-type benchmark using strong-form VSL, where the optimized objective is the quadrature approximation of the residual energy $E_{\text{strong}}$, and the solution is again parameterized by a Dirichlet-satisfying Chebyshev basis $N_x = 16$. In this setting, the training dynamics showed an overall reduction of the strong residual energy from $E \approx 1.721$ initially to $E \approx 5.234 \times 10^{-3}$ by epoch 2900. The diagnostic PDE residual loss likewise decreased but exhibited pronounced oscillations, particularly after cosine-restart events (e.g., a sharp increase in PDE loss near epoch 1600 and subsequent fluctuations through later cycles). Such behavior is typical for warm restarts. In terms of solution accuracy, strong-form VSL achieved $L^2_{\text{rel}} = 5.215 \times 10^{-4}$ and $L^\infty_{\text{rel}} = 6.596 \times 10^{-4}$. The strong-form objective produced a noticeably less accurate final solution than the GLS weak-form VSL case (roughly an order of magnitude worse in $L^2_{\text{rel}}$). This outcome is consistent with the observation that pointwise (strong) least-squares residual minimization can over-emphasize localized residual structure and become more sensitive to optimization hyperparameters, restart timing, and the distribution of quadrature/collocation points—whereas the GLS weak form effectively "filters" the residual through projections onto the basis, which can be more stable for smooth solutions and moderate mode counts.

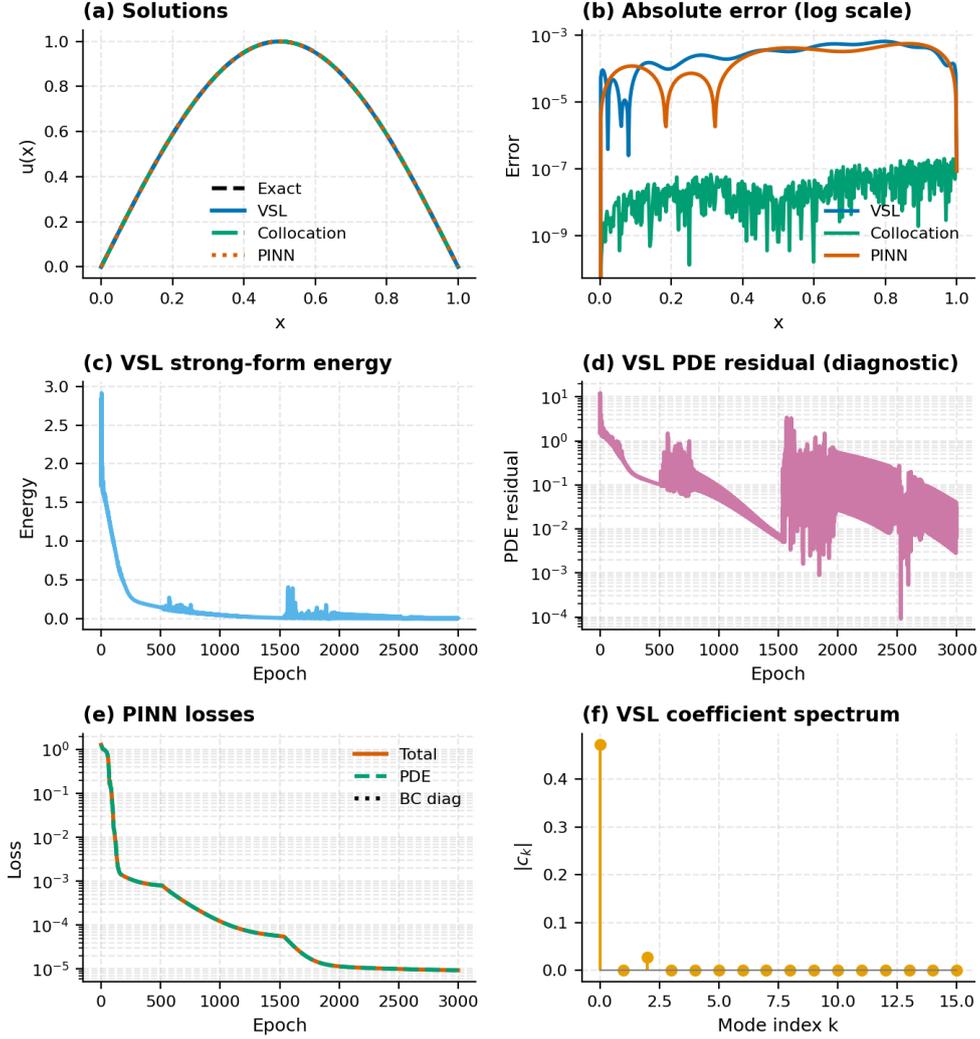

**Figure 18.** The same as in figure 17, but for VSL strong-form objective.

The two baselines behaved similarly to the weak-case campaign. classical Chebyshev collocation again delivered near machine-precision accuracy $L^2_{\text{rel}} = 7.653 \times 10^{-8}$, $L^\infty_{\text{rel}} = 2.038 \times 10^{-7}$, reinforcing that for this smooth 1D manufactured problem, the traditional spectral Newton solve remains the reference method. The PINN achieved $L^2_{\text{rel}} = 4.466 \times 10^{-4}$, $L^\infty_{\text{rel}} = 5.678 \times 10^{-4}$, with boundary lifting maintaining negligible BC error. Interestingly, in this run the PINN slightly outperformed strong-form VSL in $L^2_{\text{rel}}$, highlighting that strong-form coefficient learning is not automatically superior and must be tuned carefully.

We benchmarked three solvers on the manufactured 2D steady Burgers-type boundary value problem on $\Omega = (0,1)^2$ with homogeneous Dirichlet data and viscosity $\nu = 0.1$, using the exact field $u^\star(x, y) = \sin(\pi x) \sin(\pi y)$ to define the forcing term. The comparison includes: (i) VSL in coefficient space using a Dirichlet-satisfying Chebyshev tensor basis with $N_x = N_y = 8$, trained either with a GLS weak energy or a strong (pointwise residual) energy; (ii) a classical Chebyshev–Lobatto collocation discretization with ($N_x = N_y = 24$) solved by damped Newton; and (iii) a boundary-lifted PINN with a ($64 \times 4$) tanh MLP. Errors are reported on a uniform ($64 \times 64$) test grid. Figures 19–22 visualize the final fields (solutions and absolute errors), a midline slice at ($y \approx 0.5$), and the training diagnostics.

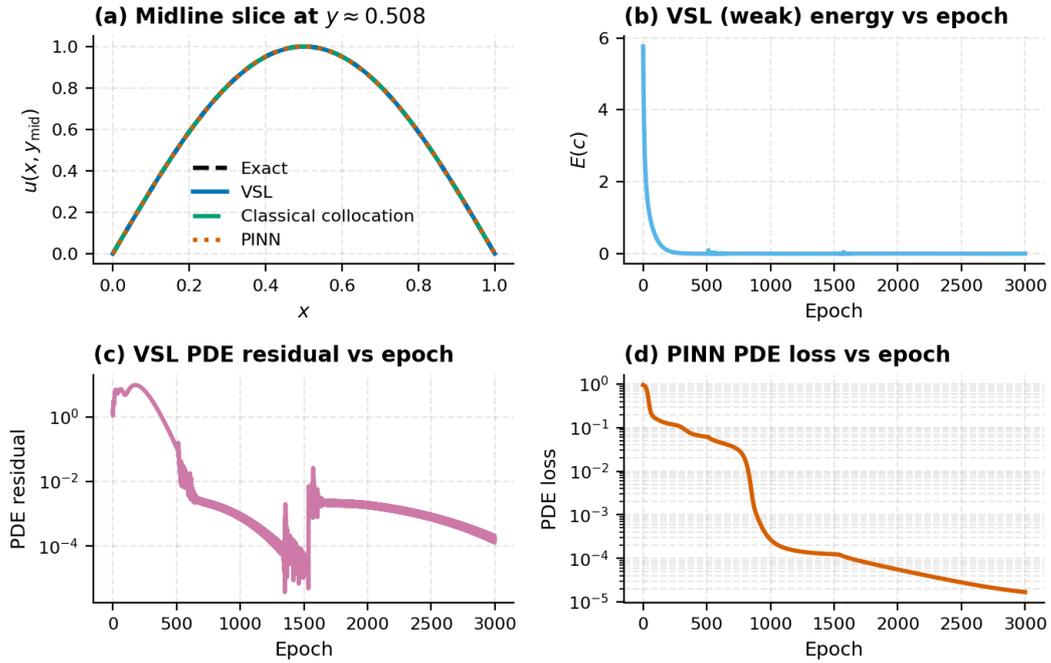

**Figure 19.** VSL (GLS weak energy). Midline validation and training diagnostics $\nu = 0.1$. (a) Midline slice $u(x, y_{\mathrm{mid}})$ at $y_{\mathrm{mid}} \approx 0.5$, comparing VSL, classical collocation, and PINN against the exact solution $u^\star$. (b) VSL training history of the GLS weak energy. (c) VSL mean-square strong residual evaluated on interior diagnostic collocation points (log scale). (d) PINN PDE loss over the same collocation set (log scale).

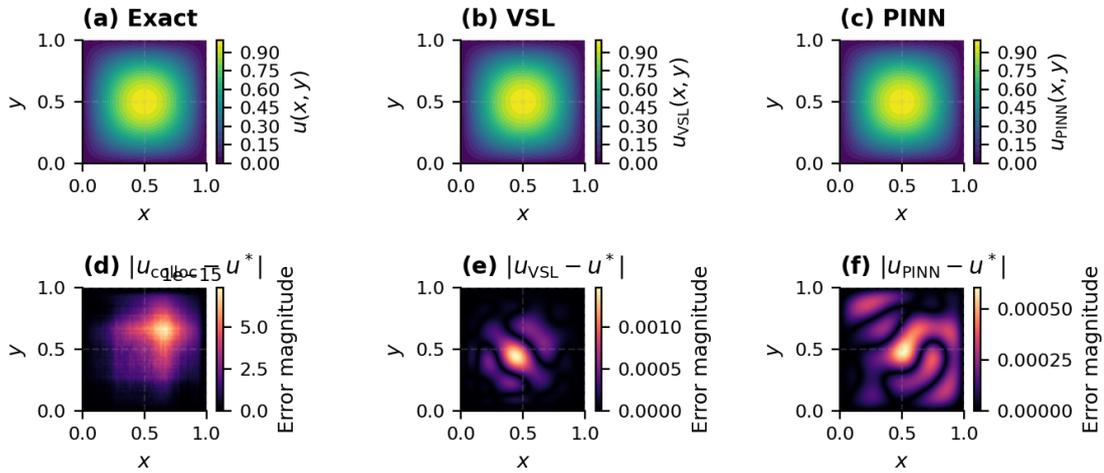

**Figure 20.** VSL (GLS weak energy). Solutions and error fields for the 2D steady Burgers-type problem ($\nu = 0.1$). Top row: filled-contour plots of (a) the manufactured exact solution $u^\star(x, y) = \sin(\pi x) \sin(\pi y)$, (b) the VSL solution $u_{\mathrm{VSL}}$ obtained by minimizing the GLS weak energy, and (c) the boundary-lifted PINN solution $u_{\mathrm{PINN}}$. Bottom row: absolute error fields on the uniform test grid—(d) classical Chebyshev–Lobatto collocation error $|u_{\mathrm{colloc}} - u^\star|$, (e) VSL error $|u_{\mathrm{VSL}} - u^\star|$, and (f) PINN error $|u_{\mathrm{PINN}} - u^\star|$.

## Midline slice and training diagnostics (VSL & PINN)

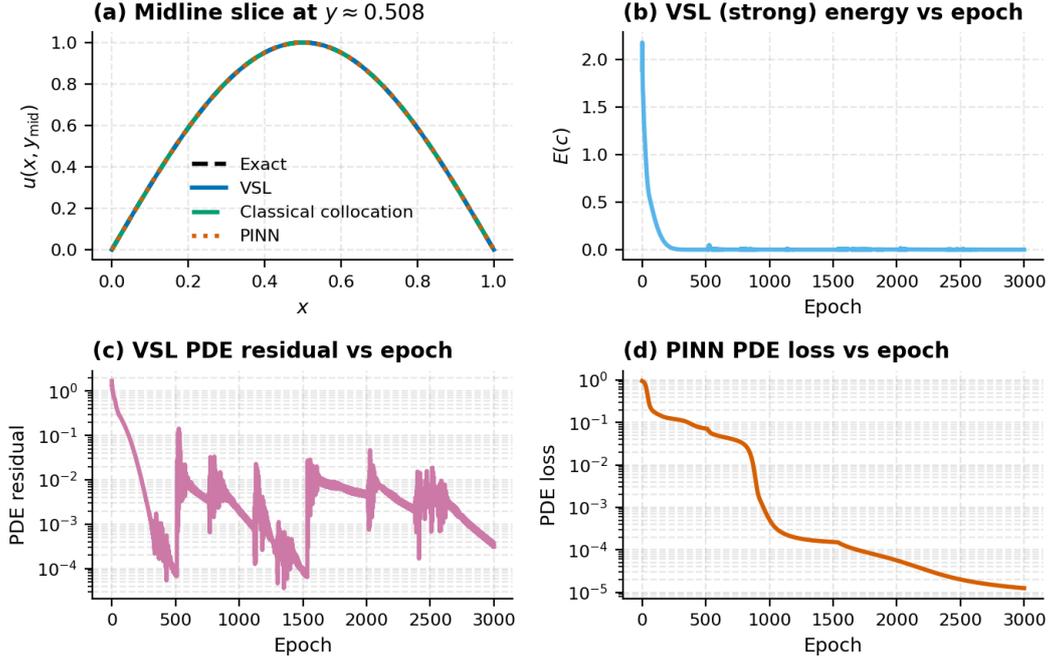

**Figure 21.** The same as in figure 19, but for VSL strong-form objective.

## Solutions and error fields

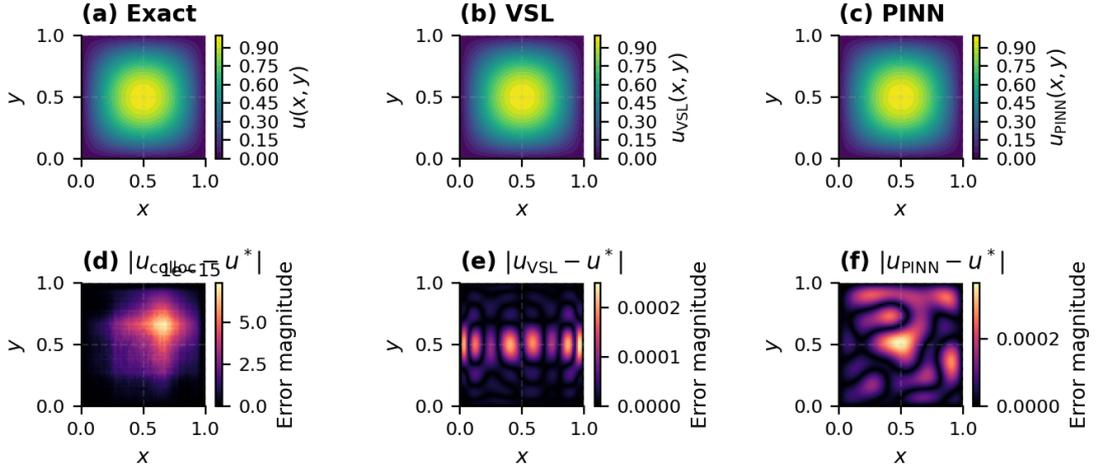

**Figure 22.** The same as in figure 20, but for VSL strong-form objective.

In the GLS Weak configuration, VSL minimizes the squared norm of the projected (Galerkin) residual vector, rather than penalizing the pointwise residual directly. The training history shows a rapid initial drop of the weak energy $E_{\text{weak(GLS)}}$ by several orders of magnitude, followed by non-monotone behavior associated with cosine-decay restarts (visible as jumps around epochs $\sim 600$ and $\sim 1600$). The diagnostic PDE residual measured at interior collocation points decreases to the $10^{-4} - 10^{-3}$ range by the end of training. On the test grid, GLS weak VSL achieves $L^2_{\text{rel}} = 5.696 \times 10^{-4}$ and $L^\infty_{\text{rel}} = 1.470 \times 10^{-3}$.

In the strong configuration, VSL minimizes the quadrature approximation of the residual least-squares energy. This objective aligns more directly with pointwise PDE satisfaction and yields noticeably improved accuracy with the same spectral model size. The training curves again show restart-induced non-monotonicity, but the method reaches a significantly smaller final solution error on the test grid: $L^2_{\text{rel}} = 1.166 \times 10^{-4}$ and $L^\infty_{\text{rel}} = 2.513 \times 10^{-4}$.

The classical Chebyshev–Lobatto collocation method attains near machine-precision agreement with the manufactured solution on the test grid: $L^2_{\text{rel}} \approx 4.425 \times 10^{-15}$ and $L^\infty_{\text{rel}} \approx 7.332 \times 10^{-15}$. This is expected for a smooth manufactured solution combined with a high-order spectral discretization and a Newton solve of the resulting nonlinear algebraic system. The PINN converges steadily in PDE loss and delivers competitive accuracy among learned methods. In the weak-case run, PINN reaches $L_{\text{rel}}^2 \approx 3.676 \times 10^{-4}$ and $L_{\text{rel}}^\infty \approx 6.056 \times 10^{-4}$, while in the strong-case run it improves to $L_{\text{rel}}^2 \approx 2.364 \times 10^{-4}$ and $L_{\text{rel}}^\infty \approx 3.670 \times 10^{-4}$.

## 9. Summary of contributions, Key lessons, and Future work

The work presented in this article has introduced and systematically evaluated VSL as a machine-learning framework for variational PDE solvers, grounded explicitly in classical functional-analytic formulations while leveraging modern optimization and automatic differentiation tools. At its core, VSL returns to the variational roots of PDE discretization by constructing space–time energies derived from either strong-form least-squares formulations or Galerkin-type weak residual formulations. In the latter case, "weak" residuals are implemented either through integration-by-parts weak forms when appropriate (e.g., for coercive elliptic/parabolic operators) or through moment residual projections of the strong residual against the trial basis, depending on the benchmark and operator structure. The unknown solution $u$ is parameterized in a tensor-product spectral basis that analytically enforces homogeneous Dirichlet boundary conditions in space: $u_N$ is represented as a linear combination of boundary-satisfying spatial basis functions and temporal Chebyshev polynomials, and the corresponding spectral coefficients $\mathbf{c}$ are treated as trainable parameters. The central objective takes the form $\mathcal{J}(\mathbf{c}) = E_{\text{strong/weak}}(\mathbf{c}) + \lambda_{\text{IC}} \mathcal{L}_{\text{IC}}(\mathbf{c}) + \lambda_{\text{reg}} \|\mathbf{c}\|_2^2$ where $E_{\text{strong/weak}}(\mathbf{c})$ encodes the PDE constraints in space–time, $\mathcal{L}_{\text{IC}}(\mathbf{c})$ enforces initial conditions when present, and $\lambda_{\text{IC}}$, $\lambda_{\text{reg}} > 0$ weight the initial-condition and regularization terms. Within this setting, VSL is implemented in a standard deep-learning framework using automatic differentiation to assemble residuals and gradients, and optimized with gradient-based methods equipped with a cosine-decay-with-restarts learning-rate schedule. The framework is evaluated across three representative benchmark families: elliptic Poisson problems in one and two dimensions, one- and two-dimensional time-dependent diffusion equations treated in a fully space–time fashion, and Burgers-type diffusion–advection problems that introduce nonlinear transport effects. In each case, VSL is compared against classical Chebyshev spectral collocation solvers and PINNs with analytic boundary–initial condition lifting, yielding a coherent experimental suite that elucidates coefficient-space learning relative to both traditional solvers and field-space machine-learning approaches.

Several methodological contributions emerge from this study. First, the construction of strong-form and Galerkin-type residual energies in coefficient space—combined with tensor-product Gauss–Legendre quadrature—provides a systematic recipe for building variational losses that remain fully compatible with automatic differentiation. The same computational machinery applies across elliptic, parabolic, and mildly nonlinear PDEs; only the differential operator and forcing definition change between benchmarks. Second, the use of Dirichlet-satisfying Chebyshev bases in space enforces homogeneous boundary conditions analytically: for any coefficient vector $\mathbf{c}$, the approximation satisfies the prescribed spatial boundary conditions by construction and lies in the intended trial space. Moreover, the polynomial basis provides sufficient smoothness to evaluate the derivatives required by strong-form energies (e.g., $u_t$, $u_{xx}$, $u_{yy}$) in a consistent pointwise or $L^2$ sense. This eliminates boundary penalty terms and mitigates the delicate loss-balancing issues often encountered in PINN training. Third, the space–time spectral parameterization for time-dependent problems removes the need for explicit time stepping and recasts evolution as a single variational minimization over coefficient space. This contrasts with classical time-marching schemes that separate spatial discretization from temporal integration and with many PINN formulations that do not explicitly control the approximation space. Finally, by comparing VSL, PINNs, and classical Chebyshev collocation (including Crank–Nicolson-type baselines), the study provides a structured reference for interpreting accuracy, conditioning, and optimization behavior under matched manufactured-solution settings where errors and diagnostics can be quantified reliably.

The numerical experiments suggest several lessons and practical recommendations for VSL-type methods. A first lesson is that VSL can match, and in favorable regimes closely approximate, the accuracy of classical high-order spectral solvers for smooth problems when the spectral resolution and quadrature rules are chosen commensurately and the optimization is driven sufficiently close to stationarity.

A second lesson is that the distinction between strong-form and weak-form residual energies matters in practice. For Poisson and diffusion problems, strong-form least-squares energies are straightforward to implement and often yield robust convergence, but they involve higher-order derivatives and can be more sensitive to quadrature accuracy and conditioning. Galerkin-type weak residual energies typically involve only first derivatives (after integration by parts) and can exhibit more favorable numerical behavior for coercive elliptic and parabolic operators, particularly when the underlying bilinear form is symmetric and coercive. The experiments indicate that both formulations are viable within VSL, and that the preferable choice depends on PDE type, regularity, derivative order, and the attainable quadrature accuracy on the chosen domain.

A third lesson concerns the role of global space–time representations for time-dependent problems. For the two-dimensional diffusion equation, the fully space–time VSL formulation circumvents explicit time stepping and constructs an approximation $u_N(x, y, t)$ that is globally valid over the space–time cylinder. This global viewpoint aligns naturally with the optimization perspective: training minimizes a single energy aggregating residual information across all times, rather than propagating information step-by-step. When the solution is smooth and the time interval is moderate, this approach achieves competitive accuracy relative to a classical Chebyshev–Lobatto spatial discretization combined with Crank–Nicolson time stepping, and the evolution of the energy and diagnostic residuals during training provides informative convergence indicators. At the same time, the experiments suggest that space–time spectral representations may require careful temporal basis selection and regularization when pronounced temporal stiffness or multiple separated time scales are present; in such regimes, a larger number of temporal modes may be needed and the optimization landscape can become more delicate. Practically, the temporal spectral resolution should therefore be treated as a tunable hyperparameter, guided by diagnostics such as energy decay, residual monitoring on independent collocation sets, and sensitivity studies.

A fourth lesson is that coefficient-space learning in VSL offers a controlled and interpretable alternative to overparameterized field-space PINNs in the settings examined. Because the VSL parameterization is linear in **c** with a fixed, analytically known feature map (the spectral basis), the nonlinearity resides entirely in the PDE residual and its dependence on **c**. This yields optimization problems that resemble deterministic nonlinear least-squares in a structured feature space, rather than the highly non-convex landscapes associated with deep neural parameterizations. Empirically, VSL exhibits more predictable training dynamics and reduced sensitivity to architectural hyperparameters such as width, depth, and activation choice, which are absent by construction. In contrast, PINN baselines often require more careful tuning to achieve comparable residual reduction, particularly in Burgers-type problems where nonlinear transport amplifies gradient pathologies. At the same time, the experiments reaffirm an important trade-off: PINNs retain greater flexibility in principle for complex geometries and localized features without explicit basis design, whereas VSL performance is closely tied to the suitability of the chosen spectral space. This suggests that VSL is especially effective as a high-order, interpretable solver on simple or moderately complex domains and as a strong baseline for benchmarking, while PINNs and other field-space models remain attractive in settings where geometry, coefficients, or localized features preclude simple global bases.

A fifth lesson pertains to optimization and learning-rate schedules. By implementing VSL within a standard ML framework and using cosine-decay-with-restarts schedules, the study shows that optimization practices developed for deep learning can be profitably repurposed for variational PDE solvers in coefficient space. Such schedules can improve early-stage exploration of the objective and reduce stagnation in plateau-like regions, while still providing progressive step-size reduction for late-stage refinement. Across the benchmarks, these schedules support stable convergence for both VSL and PINN models, with VSL typically converging in fewer effective epochs due to its reduced parameter count and smoother landscape. Nevertheless, the Burgers-type experiments underscore that nonlinearity can introduce slow-decay phases and increased sensitivity to schedule parameters, particularly when strong-form energies are used. A practical recommendation is therefore to couple VSL with robust diagnostic monitoring of energy, residuals, and learning-rate trajectories, and to treat schedule parameters (initial learning rate, restart frequency, and minimum learning-rate fraction) as meaningful hyperparameters that may require modest tuning across PDE families.

The findings of this work open several avenues for future research. A natural extension is to move beyond rectangular domains and homogeneous coefficients toward more complex geometries and heterogeneous media. Achieving this while preserving the "boundary-by-construction" philosophy will likely require mapped coordinates, patchwise local bases, or constrained trace spaces, potentially combined with domain decomposition. Incorporating curvilinear mappings, spectral elements, or hybrid discretizations that combine local spectral bases with finite element or finite volume components could broaden applicability while retaining the coefficient-space learning paradigm. Relatedly, parametric PDEs and uncertainty quantification [69-71] represent an important application area: one may envisage learning mappings from parameter spaces to coefficient vectors, thereby integrating VSL as a variationally grounded building block within operator-learning architectures.

Another promising direction concerns adaptivity and basis enrichment. The present study employs fixed tensor-product Chebyshev bases, which are highly effective for smooth solutions but can become inefficient when localized structures, boundary layers, or near-discontinuities arise. Extending VSL to incorporate adaptive basis selection—through hierarchical spectral families, multi-resolution representations, or localized enrichment driven by residual indicators—could substantially broaden the method's scope. This adaptivity introduces algorithmic questions about how to expand the parameter space during training while preserving numerical stability, and it raises theoretical questions about the interplay between basis enrichment, coefficient regularization, and learning-rate scheduling.

From a theoretical standpoint, there is scope for deeper analysis of VSL convergence when the energy $J(\mathbf{c})$ is minimized by gradient-based optimization rather than solved via direct or iterative linear algebra. For linear elliptic and parabolic problems, the strong- and weak-form energies are convex quadratics under standard assumptions, and classical theory guarantees convergence of deterministic gradient descent under suitable step-size conditions. In contrast, the use of adaptive optimizers, stochastic mini-batching, additional loss terms, or nonlinear operators (as in Burgers-type problems) motivates a more refined analysis of practical convergence behavior and stopping criteria. Establishing convergence-rate bounds under realistic optimization regimes, quantifying the effect of learning-rate schedules on error, and relating coefficient-space iterates to classical Galerkin solutions are all directions that merit rigorous study.

Finally, VSL is best viewed as a complementary point in the design space of PDE solvers—bridging classical variational discretization and modern ML tooling—rather than as a universal replacement for established numerical schemes or PINN frameworks. Its strengths lie in its tight integration with variational theory, its high-order spectral approximation spaces, and its compatibility with ML software stacks that provide automatic differentiation, robust optimizers, and scheduling utilities. As such, VSL provides a useful reference point for future developments in AI for scientific computing: it offers a concrete and reproducible recipe for constructing variational, coefficient-space solvers; it supports fair comparisons with field-space PINNs and classical spectral schemes on common benchmark families; and it illustrates how numerical analysis and deep learning concepts can be combined in a principled manner. Further work extending VSL to broader PDE classes, refining its numerical analysis, and exploring hybrid architectures that blend spectral representations with learned components has the potential to deepen our understanding of how machine learning can contribute to the design of accurate, robust, and interpretable PDE solvers.

## Code availability

All code used to generate the data, figures, and analyses in this study is publicly available, together with implementation details, at the Zenodo repository (https://doi.org/10.5281/zenodo.18156105) [72].